\documentclass[11pt]{amsart}
\usepackage{amsmath}
\usepackage{amssymb}
\usepackage{mathrsfs}

\newcommand{\conc}{{}^\frown\!}
\newcommand{\lh}{\ell g\/} 
\newcommand{\rest}{{\restriction}}
 
\newcommand{\rng}{{\rm rng}}
\newcommand{\cf}{{\rm cf}}
\newcommand{\bt}{{\bf t}}
\newcommand{\D}{{\mathcal D}}

\newcommand{\Hla}{{\mathcal H}(\lambda)}
\newcommand{\DOM}{{\rm DOM}}
\newcommand{\ID}{{\rm ID}}
\newcommand{\id}{{\rm id}}
\newcommand{\wdmid}{{{\rm WDmId}_\lambda}}
\newcommand{\wdmidpl}{{{\rm WDmId}_{\lambda^+}}}
\newcommand{\MOD}{{\rm MOD}}
\newcommand{\Frag}{{\rm Frag}}
\newcommand{\gK}{{\mathfrak K}}
\newcommand{\gB}{{\mathfrak B}}
\newcommand{\ideal}{{\rm ideal}_\lambda}
\newcommand{\dun}{\mathop{\nabla}}
\newcommand{\tr}{{\rm tr}}
\newcommand{\lqq}{{\text{`}\text{`}}} 
\newcommand{\ida}{{\rm ID^a}}
\newcommand{\idb}{{\rm ID^b}}
\newcommand{\idc}{{\rm ID^c}}
 \newcommand{\idi}{{{\rm ID}^ \iota }}

\newtheorem{theorem}{Theorem}[section] 
 
\newtheorem{claim}{Claim}[theorem]
\newtheorem{proposition}[theorem]{Proposition}

\theoremstyle{definition}
\newtheorem{definition}[theorem]{Definition}

\newtheorem{hypothesis}[theorem]{Hypothesis}
 
\newtheorem{construction}[theorem]{Construction}

\theoremstyle{remark}
\newtheorem{conclusion}[theorem]{Conclusion}
\newtheorem{remark}[theorem]{Remark}
\newtheorem{notation}[theorem]{Notation}

\title{More on Weak Diamond}
\author{Saharon Shelah}
\address{Institute of Mathematics\\
 The Hebrew University of Jerusalem\\
 91904 Jerusalem, Israel\\
 and  Department of Mathematics\\
 Rutgers University\\
 New Brunswick, NJ 08854, USA\\
 and  Mathematics Department\\
 University of Wisconsin -- Madison\\
 Madison, WI 53706, USA}
\email{shelah@math.huji.ac.il}
\urladdr{http://www.math.rutgers.edu/$\sim$shelah}
\thanks{The research partially supported by The National Science
Foundation Grant No. 144--EF67. Publication 638}   

\date{2021-03-03a}

\usepackage[hidelinks]{hyperref}
\begin{document}
\makeatletter\def\shfiuwefootnote{\gdef\@thefnmark{}\@footnotetext}\makeatother\shfiuwefootnote{Version 2021-03-03\_2. See \url{https://shelah.logic.at/papers/638/} for possible updates.} 
\begin{abstract}
We deal with the combinatorial principle Weak Diamond.
We prove that, if it holds 
for a given cardinal, 
we can get this principle with more than two colours
or some relevant ideal is not too saturated.
Then we point out a model theoretic consequence of Weak Diamond.
\end{abstract}
\maketitle 

\section{Basic definitions}
In this section we present basic notations, definitions and results. 

The paper was circulated (including the math arXiv) and accepted to the 
East-West Journal of Math around 2000, but due to some problems
between the editors has not appeared.  Meanwhile Aspero, Larson and Moore
\cite{AsLaMo13} 
with a related result was  done. 
Weak diamond was introduced in \cite{Sh:65}, lately see
\cite{Sh:1111}.

\begin{notation}  \label{z2}
\begin{enumerate}
\item $\kappa,\lambda,\theta,\mu$ will denote cardinal numbers and $\alpha$,
$\beta$, $\delta$, $\varepsilon$, $\xi$, $\zeta$, $\gamma$ will be
used to denote ordinals.  
\item Sequences of ordinals are denoted by $\nu$, $\eta$, $\rho$ (with
possible indexes).
\item The length of a sequence $\eta$ is $\lh(\eta)$.
\item For a sequence $\eta$ and $\ell\leq\lh(\eta)$, $\eta\rest\ell$ is the
restriction of the sequence $\eta$ to $\ell$ (so $\lh(\eta\rest\ell)=\ell$). 
If a sequence $\nu$ is a proper initial segment of a sequence $\eta$ then we
write $\nu\vartriangleleft\eta$ (and $\nu\trianglelefteq\eta$ has the obvious
meaning).  
\item For a set $A$ and an ordinal $\alpha$, $\alpha_A$ stands for the
function on $A$ which is constantly equal to $\alpha$. 
\item For a model $M$, $|M|$ stands for the universe of the model.
\item The cardinality of a set $X$ is denoted by $\|X\|$. The cardinality of
the universe of a model $M$ is denoted by $\|M\|$.
\end{enumerate}
\end{notation}

\begin{definition}  \label{z5}  
\label{colourings}
Let $\lambda$ be a regular uncountable cardinal and $\theta$ be a cardinal
number, possibly finite.  
\begin{enumerate}
\item {\em A $(\lambda,\theta)$--colouring} is a function $F:\DOM
\longrightarrow\theta$, where $\DOM$ is either ${}^{\textstyle {<}\lambda}2=
\bigcup\limits_{\alpha<\lambda} {}^{\textstyle\alpha}2$ or $\bigcup\limits_{
\alpha<\lambda}{}^{\textstyle \alpha}({\mathcal H}(\lambda))$. In the first
case we will write $\DOM_\alpha=  \DOM_ \alpha (F) =  
{}^{\textstyle 1{+}\alpha}2$, in the second
case we let $\DOM_\alpha=  \DOM_\alpha ( F)  = 
{}^{\textstyle 1{+}\alpha}({\mathcal H}(\lambda))$
(for $\alpha\leq\lambda$).

If the choice does not matter we shall not mention it; for the main 
definitions the choice does not matter, see \ref{a26}.  

If $\lambda$ is understood 
from the context,  
we may omit it; if $\theta=2$ then we may omit
it 
(thus 
a $ \lambda  $-colouring means  
a $ ( \lambda, 2 ) $-colouring and 
{\em a colouring} is a $(\lambda,2)$--colouring).

\item For a $(\lambda,\theta)$--colouring $F$ and a set $S\subseteq\lambda$,
we say that a function $\eta\in {}^{\textstyle S} \theta$ is {\em an $F$--weak
diamond sequence for $S$} \underline{when} 
for every $f\in\DOM_\lambda$ the set
\[\{\delta\in S:\eta(\delta)=F(f\rest\delta)\}\]
is stationary.

\item $\wdmid$ is the collection of all sets $S\subseteq\lambda$ such that for
some colouring $F$ there is no $F$--weak diamond sequence for $S$. 
\end{enumerate}
\end{definition}

\begin{remark}   \label{z8}
In the definition of $\wdmid$ (\ref{colourings}(3)), the choice of $\DOM$
(see \ref{colourings}(1)) does not matter; see \cite[AP, \S1]{Sh:f},
remember that $\|\Hla\|=2^{{<}\lambda}$.
\end{remark}

\begin{theorem}  \label{z11}
[Devlin and   
Shelah \cite{Sh:65}; see {\cite[AP, \S1]{Sh:f}} too]
\ \\ 
Assume that $2^\theta=2^{<\lambda}<2^\lambda$ (e.g.\ $\lambda=\mu^+$,
$2^\mu< 2^\lambda$). Then for every $\lambda$--colouring $F$ there exists an
$F$--weak diamond sequence for $\lambda$. Moreover, $\wdmid$ is a normal
ideal on $\lambda$ (and $\lambda\notin\wdmid$).
\end{theorem}

\begin{remark}  \label{z14}
One could wonder why the weak diamond (and $\wdmid$) is interesting. Below we
list some of the applications, limitations and related problems.
\begin{enumerate}
\item Weak diamond is really weaker than diamond, but 
provably (in ZFC)  
it holds true for some
cardinals $\lambda$.
Note that under GCH, $\diamondsuit_{\mu^+}$ holds
true for each $\mu>\aleph_0$, so the only interesting case then is
$\lambda=\aleph_1$. 
\item Original interest in this combinatorial principle comes from
interest in   
Whitehead
groups: 
\begin{quotation}
{\em \noindent if $G$ is a strongly $\lambda$--free Abelian group and
$\Gamma(G)\notin \wdmid$

\noindent then $G$ is  not  Whitehead.}
\end{quotation}
\item A related question was:\quad can we have stationary subsets $S_1,S_2
\subseteq\omega_1$ such that $\diamondsuit_{S_1}$ but
$\neg\diamondsuit_{S_2}$? (See \cite{Sh:64}.)
\item Weak diamond has been helpful particularly in problems where we have
some uniformity, e.g.:
\begin{enumerate}
\item[$(*)_1$]{\em Assume $2^\lambda<2^{\lambda^+}$. Let
 $\psi\in {\mathbb   L}   _ { 
\lambda^+,\omega}$ be categorical in $\lambda,\lambda^+$.\\
Then $(\MOD_\psi,\prec_{\Frag(\psi)})$ has the amalgamation property in
$\lambda$.} 
\item[$(*)_2$]{\em If $G$ is a 
group  of cardinality $ \lambda >   
{\aleph_0} $ 
then we can find subgroups
$G_i$ of $G$ (for $i<\lambda$) non-conjugate in pairs (see \cite{Sh:192}).}
\end{enumerate}
\item One may wonder if assuming $\lambda=\mu^+$, $2^\lambda>2^\mu$ (and
e.g.\ $\mu$ regular) we may find a regular $\sigma<\mu$ such that 
\[
\{\delta<\lambda:\cf(\delta)=\sigma\}\notin\wdmid(\lambda).
\] 
Of course, by the \lqq normal ideal"    
   result  
   (see   
\ref{z11}), it follows that there is such $ \sigma $, but does 
$ \sigma $ depends   
on the present set theory? 
e.g. does it hold   
for every regular $ \sigma \not= \cf(\mu )$ 
below $ \mu $?   

Unfortunately, this is not the case (see \cite{Sh:208}
even for $ \mu = {\aleph_1} $).
\item We would like to prove
\begin{enumerate}
\item[(a)] $\wdmid$ is not $\lambda^+$--saturated\quad or
\item[(b)] a strengthening, e.g.\ weak diamond for more 
(than two) 
colours.
\end{enumerate}
We will get (a variant of) a local version of the disjunction, where we
essentially fix $F$. There are two reasons for interest in {\bf (a)}:
understanding $\lambda^+$--saturated normal ideals (e.g.\ we get more
information on the case CH + ``$\D_{\omega_1}$ is $\aleph_2$--saturated''; see
also Zapletal Shelah \cite{Sh:610}), and non $\lambda^+$--saturation helps
in ``non-structure theorems in model thery''   
(see \cite{Sh:87b}, \cite{Sh:576}, 
\cite{Sh:h}, \cite{Sh:i}).  
That is,
having $2^\mu<2^{\mu^+}<2^{\mu^{++}}$ and some ``bad'' (i.e.\
``non-structure'') properties for models in $\mu$ we get $2^{\mu^{++}}$ models
in $\mu^{++}$ when $\wdmidpl$ is not $\lambda^{++}$--saturated (and using the
local version does not hurt).
\item Note that for $S\notin\wdmid$ we have a weak diamond sequence $f\in
{}^{\textstyle S}2$ such that the set of ``successes''(=equalities) is
stationary, but it does not have to be in $(\wdmid)^+$. We would like to
start and end in the same place: being positive for the same ideal. 
Also, in {\bf (b)} above the set of places we guess was stationary, when
we start with $S\in (\wdmid)^+$. 

Note that it may well be that $\lambda\in\wdmid$ (if $(\exists\theta<
\lambda)(2^\theta=2^\lambda)$ this holds), but some ``local'' versions may
still hold. E.g.\ in the Easton model, we have $F$--weak diamond sequences
for all $F$ which are reasonably definable (see \cite[AP, \S1]{Sh:f}; 
define   e.g.  
\[F(f)=1\ \Leftrightarrow\ {\rm L}[X,f]\models\varphi(X,f)\]
for a fixed first order formula $\varphi$, where $X\subseteq\lambda$ depends
on $F$ only). So the case $\wdmid={\mathcal P}(\lambda)$ has some interest. 
\item
Related later works  are  \cite{Sh:897},  \cite{Sh:1028}.   
\end{enumerate}
\end{remark}

We would like to thank Andrzej Ros{\l}anowski for mathematical comments and
improving the presentation 
and to Shimoni Garti.  

\section{When colourings are almost constant}

In \ref{a2}, \ref{a5} we now \lqq slice"  $ \wdmid $ 
by finer approximations meaningful even when $ \wdmid $  fails.  

\begin{definition}   \label{a2}
Let $\lambda$ be a regular uncountable cardinal.
\begin{enumerate}
\item Let    
$F$ be a $(\lambda,\theta)$--colouring. 
\begin{enumerate} 
    \item[(a)] Let  $S\subseteq\lambda$ . 
We say that {\em a sequence $\eta\in {}^{\textstyle S}\theta$ is coded by $F$}
if there exists $f\in\DOM_\lambda$ such that 
\[\alpha\in S\quad\Leftrightarrow\quad \eta(\alpha)= F(f\rest (1+\alpha)).\]
    \item[(b)] We let 
\[\gB(F)\stackrel{\rm def}{=}\{\eta\in {}^{\textstyle\lambda}\theta:
\eta\mbox{ is coded by }F\/\}.\]
\end{enumerate} 
\item For a family ${\mathcal A}$ of subsets of $\lambda$ let
$\ideal({\mathcal A})$ be the $\lambda$--complete normal ideal on $\lambda$
generated by ${\mathcal A}$ (i.e.\ it is the closure of ${\mathcal A}$ under
unions of $<\lambda$ elements, diagonal unions, containing singletons, and
subsets). 

[Note that $\ideal({\mathcal A})$ does not have to be a proper ideal.]
\item For a $\lambda$--colouring $F$ (so $\theta=2$) we define by induction on
$\alpha$: 
\[\ID^-_0(F)=\emptyset,\qquad \ID_0(F)=\{S\subseteq\lambda: S\mbox{ is not
stationary }\},\] 
for a limit $\alpha$
\[\ID^-_\alpha(F)=\bigcup_{\beta<\alpha}\ID_\beta(F),\qquad \ID_\alpha(F)
=\ideal(\bigcup\limits_{\beta<\alpha}\ID_\beta(F)),\] 
and\footnote{ 
    Note that 
    $ \ID^-_\alpha (F)  \not= \emptyset $  iff $ \alpha > 0 $  
    }
for $\alpha=\beta+1$
\[\begin{array}{ll}
\ID^-_\alpha(F)=&\big\{S\subseteq\lambda:\mbox{for each }S^*\subseteq S\mbox{
there is }f\in \DOM_\lambda\mbox{ such that}\\
\ &\qquad\qquad\{\delta<\lambda: \delta\in S^*\Leftrightarrow F(f\rest\delta)
=0\}\in\ID_\beta(F)\big\};\\
\ID_\alpha(F)=&\ideal\big(\ID^-_\alpha(F)).
  \end{array}\]
Finally we let $\ida(F)=\bigcup\limits_{\alpha}\ID_\alpha(F)$.
\item We say that $F$ is {\em rich\/} \underline{when}   
if 
$\DOM(F)=\bigcup\limits_{\alpha<
\lambda}{}^{\textstyle \alpha}{\mathcal H}(\lambda)$, and for every function
$f\in\DOM_\lambda$ and $\alpha<\lambda$ and a set $A\subseteq\alpha$ there is
$f'\in\DOM_\lambda$ such that:
\[
(\forall i   \in ( \lambda \setminus  \alpha )
    (f(1+i)=f'(1+i)\ \&\ F(f\rest (\alpha+i))=F(f'\rest
(\alpha+i)))
\]
and $(\forall j   \in \alpha \setminus \{ 0\}  
)(F(f'\rest j)=1\ \Leftrightarrow\ j\in A)$.
\end{enumerate}
\end{definition}

\begin{definition}  \label{a5}
\label{1.1Adef}
Let $\lambda$ be a regular uncountable cardinal and let $F$ be a
$\lambda$--colouring.  
\begin{enumerate}
    \item $\wdmid(F)$ is the family of all sets $S\subseteq\lambda$ with the
property that for every $S^*\subseteq S$ there is $f\in\DOM_\lambda$ such that
the set
\[\{\delta\in S:\delta\in S^*\ \Leftrightarrow\ F(f\rest\delta)=1\}\]
is not stationary, 
    (note,  the difference with \ref{z5}(3)).  
    \item  
$\gB^+(F)$ is the closure of 
   (see \ref{a2}(1)):  
\[\gB(F)\cup\{S\subseteq\lambda: S\mbox{ is not stationary }\}\]
under unions of $<\lambda$ sets, complement and diagonal unions (here, in
$\gB(F)$, we identify a subset of $\lambda$ with its characteristic
function). 
    \item $ \idb \stackrel{\rm def}{=}\{  S \subseteq \lambda : $  
  for every $ Y \subseteq S $ for some $ B \in \mathfrak{B} ^+(F) $
  we have $ B \cap S = Y  \} $
   \item $\idc(F)$ is the collection of all $S\subseteq\lambda$ such that for
some $X\in\gB^+(F)$ we have:\quad $S\subseteq X$ and there is a partition
$X_0,X_1$ of $X$ such that 
\begin{enumerate}
\item[$(\alpha)$] ${\mathcal P}(X_\ell)=\{Y\cap X_\ell: Y\in \gB^+(F)\}$ for
$\ell=0,1$, and 
\item[$(\beta)$] for $\ell<2$   
there is no $Y\in\gB^+(F)$  satisfying
\[Y\setminus X_\ell\in\idb(F)\quad\&\quad Y\notin\idb(F).\]
\end{enumerate}
\end{enumerate}
\end{definition}

\begin{proposition}  \label{a8}
\label{easyprop}
Assume $\lambda$ is a regular uncountable cardinal and $F$ is a
$\lambda$--colouring. 
\begin{enumerate}
    \item 
$\ideal ({\mathcal A})$ is the collection of all diagonal unions
$\dun\limits_{\xi<\lambda} A_\xi$ such that $A_\xi\in{\mathcal A}$ for
$\xi<\lambda$,  
\underline{when}  
 ${\mathcal A}$ is a family of subsets of $\lambda$ such that
  \begin{enumerate}
  \item[$(\circledast_{\mathcal A})$] \quad if $S_0\subseteq S_1$ and $S_1
  \in {\mathcal A}$ and $A\in [\lambda]^{\textstyle <\!\lambda}$ then
  $S_0\cup A\in{\mathcal A}$, 
  \end{enumerate}
    \item The condition $(\circledast_{\ID^-_\alpha(F)})$ (see above) 
is true
for each $\alpha$. Consequently, if $\alpha=\beta+1$ then $\ID_\alpha(F)=\{
\dun\limits_{i<\lambda}A_i:\langle A_i:i<\lambda\rangle\subseteq\ID^-_\alpha
(F)\}$, and if $\alpha$ is limit then $\ID_\alpha(F)=\{\dun\limits_{i<
\lambda} A_i:\langle A_i:i<\lambda\rangle\subseteq\bigcup\limits_{\beta<
\alpha}\ID_\beta(F)\}$.  
    \item $\ida(F)$ and $\ID_\alpha(F)$ are $\lambda$--complete normal ideals on
$\lambda$ extending the ideal of non-stationary subsets of $\lambda$ (but
they do not have to be proper). For $\alpha<\gamma$ we have $\ID_\alpha(F)
\subseteq\ID_\gamma(F)$ and hence $\ida(F)=\ID_\alpha(F)$ for every large
enough $\alpha<(2^\lambda)^+$.  
    \item Suppose $\bar{B}=\langle B_\ell:\ell\leq m\rangle$, where $B_\ell
\subseteq B_{\ell+1}$ (for $\ell<m$) and $B_m\in\ida(F)$. Then $\bar{B}$ has
an $F$--representation, which means that there are a well founded tree $T
\subseteq {}^{\textstyle\omega\!>}\lambda$, sequences $\langle B^\ell_\eta: 
\eta\in T,\ \ell\leq\ell_\eta\rangle$, and $\langle f^k_\eta: \eta\in T,\
k\leq k_\eta\rangle$ such that 
$k_\eta\leq\ell_\eta+1
\le m + 1    
$ and
\begin{enumerate}
 \item[(a)] $B^\ell_{\langle\rangle}=B
        _{\ell}     
$, $\ell_{\langle\rangle}=m$,
$B^\ell_\eta\subseteq B^{\ell+1}_\eta\subseteq\lambda$, $f^\ell_\eta\in
{}^{\textstyle\lambda} 2$, 
 \item[(b)] $(\forall\eta\in T\setminus\max(T))(\forall i<\lambda)(\eta\conc
\langle i\rangle\in T)$,
 \item[(c)] for each $\eta\in T\setminus\max(T)$ there is $\alpha_\eta<
\lambda$ such that for all 
$ {\ell} \le {\ell} _ \eta $, 
$\delta\in\lambda\setminus\alpha_\eta$  
we have:  
\begin{enumerate}
    \item[$(\oplus)$] $\delta\in B^\ell_\eta$\quad iff

\qquad $(\exists i<\delta)(\delta\in B^\ell_{\eta\conc\langle i\rangle})$ or

\qquad $F(f^\ell_\eta\rest\delta)=1\ \&\ \neg(\exists i<\delta)(\exists k)
(\delta\in B^k_{\eta\conc\langle i\rangle})$,
\end{enumerate}
 \item[(d)] for each $\eta\in \max(T)$, $B_\eta$ is a bounded subset of
$\lambda$ with $\min(B_\eta)>
    \sup   
(\{\eta(n):n<\lh(\eta)\})$. 
\end{enumerate}
    \item If for some $f^*\in {}^{\textstyle\lambda}2$ we have $(\forall\alpha<
\lambda)(F(f^*\rest\alpha)=0)$ then in part (4) above we can demand that
$k_\eta=\ell_\eta+1$. 
    \item If $F$ is rich then in part (4) above we can add 
\begin{enumerate}
    \item[(e)] $\alpha_\eta$=0 for $\eta\in T\setminus\max(T)$ and
$B_\eta=\emptyset$ for $\eta\in\max(T)$.
\end{enumerate}
    \item $\ida(F)$ is the minimal normal ideal 
    $ D  $   
on $\lambda$ such that there is
no $S\in
D   
^+$ satisfying
\[(\forall S^*\subseteq S)(\exists A\in\gB(F))(S^*\vartriangle A\in
     D )  
.\]

    \item If $ X \in {\mathscr P} (\lambda ) \setminus \idb(F)$   
then there is $ \eta \in {}^{ \lambda } 2 $ which is a 
weak diamond even modulo $ \idb (F)$  which mean that:

for every $ f \in \DOM(F) $ we have: 
$ \{ \delta \in X: F( f \rest \delta ) = \eta (\delta )\}   
\not= \emptyset  \mod 
    \idb (F)$.

\item 
 $\ID^1(F)=  \{ S\subseteq\lambda:(\exists X\in\gB^+(F))
 (S\subseteq X\ \&\ {\mathcal P}(X)\subseteq\gB^+(F))\}$.
 \end{enumerate} 
\end{proposition}

\begin{proof} (1)  
\quad Should be clear.

\noindent 
(2)  
If $ \alpha = 0 $ then  $(\circledast_{\mathcal A})$  
holds trivially because there is no such $ S_1$.

If $ \alpha $ is a limit ordinal  then the condition holds
because for every $ S_1 \in \ID^-_\alpha (F)$
there is $ \beta < \alpha $ such that 
$ S_ 1 \in  \ID^-_\beta  (F) $ and we can use the induction 
hyp.

Lastly, if $ \alpha = \beta + 1 $ this is easy too.

\noindent 
(3)\quad  
For the first sentence $\ID^-_\alpha (F)$ is a normal ideal by 
its definition; this implies   $\ID^a(F)$   
is a normal ideal 
by the second sentence. We still have   
to prove the second sentence.  

By induction on $\gamma<\lambda$ and then by induction on
$\alpha<\gamma$ we show that $(\forall\gamma<\lambda)(\forall\alpha<\gamma)(
\ID_\alpha(F)\subseteq\ID_\gamma(F))$. If $\gamma=1$ then this follows
immediately from definitions; similarly if $\gamma$ is limit. So suppose now
that $\gamma=\gamma_0+1$ and we proceed by induction on $\alpha\leq
\gamma_0$. There are no problems 
neither   
when $\alpha=0$ nor when $\alpha$ is
limit. So suppose that $\alpha=\beta+1<\gamma$ (so $\beta<\gamma_0$). By the
inductive hypothesis we know that $\ID_\beta(F)\subseteq\ID_{\gamma_0}(F)$. 
Let $A\in\ID_{\beta+1}(F)$. By (2) there are $A_\xi\in\ID^-_{\beta+1}$ (for
$\xi<\lambda$) such that $A=\dun\limits_{\xi<\lambda}A_\xi$. Now look at the
definition of $\ID_{\beta+1}^-(F)$: since $\ID_\beta(F)\subseteq\ID_{
\gamma_0}(F)$ we see that $A_\xi\in\ID^-_{\gamma_0+1}(F)$. Hence $A\in
\ID_\gamma$.

\noindent (4)\quad By induction on $\alpha$ we show that:  if $\bar{B}=\langle
B_\ell:\ell\leq m\rangle$, where $B_\ell\subseteq B_{\ell+1}$ (for $\ell<m$)
and $B_m\in\ID_\alpha(F)$ then $\bar{B}$ has an $F$--representation.

\noindent{\sc Case 1:}\qquad $\alpha=0$.\\
Thus the set $B_m$ is not stationary and we may pick up a club $E$ of
$\lambda$ disjoint from $B_m$. Let $E=\{\alpha_\zeta:\zeta<\lambda\}$ be the
increasing enumeration. Put $T=\{\langle\rangle\}\cup\{\langle i\rangle:i<
\lambda\}$, $\alpha_{\langle\rangle}=1$, $\ell_{\langle\rangle}=\ell_{
\langle i\rangle}=m$, $B^\ell_{\langle\rangle}=B_\ell$ and $B^\ell_{\langle
i\rangle}=B_\ell\cap\alpha_{i+1}$. Now check.

\noindent{\sc Case 2:}\qquad $\alpha$ is limit.\\
It follows from (2) that $B_\ell=\dun\limits_{i<\lambda} B_{\ell,i}$
for some $B_{\ell,i}\in\bigcup\limits_{\beta<\alpha}\ID_\beta(F)$. Let
$B_{\ell,i}'$ be defined as follows: 
\smallskip

if $i=(m+1)j+t$, $\ell<t\leq m$ then $B_{\ell,i}'=\emptyset$,

if $i=(m+1)j+t$, $t\leq m$, $t\leq \ell$ then $B_{\ell,i}'=B_{\ell,i}$.
\smallskip

\noindent Then for each $i,\ell$ we may find $\langle B^{i,\ell}_\eta,
f^{i,\ell'}_\eta,\alpha^i_\eta: \eta\in T_i,\ \ell<\ell^{i,1}_\eta, \ \ell'<
\ell^{i,2}_\eta\rangle$  satisfying clauses (a)--(d) and such that $\langle
B^{\ell,i,k}_{\langle\rangle}: k\leq k^1_\eta\rangle = \langle B'_{\ell,i}:
\ell\leq m\rangle$ (by the induction   
hypothesis). Put
\[\begin{array}{l}
T=\{\langle\rangle\}\cup\{\langle i\rangle\conc\eta:\eta\in T_i\},\\
\ell_{\langle\rangle}=m,\quad \ell_{\langle\rangle}'=0,\quad\ell_{\langle
i\rangle\conc\eta}=\ell^{i,1}_\eta,\quad \ell_{\langle i\rangle\conc\eta}
=\ell^{i,2}_\eta\\
B^\ell_{\langle\rangle}=B_\ell,\quad B^\ell_{\langle i\rangle\conc\eta}=
B^{i,\ell}_\eta,\quad f^{\ell'}_{\langle i\rangle\conc\eta}=f^{i,
\ell'}_\eta,\\
\alpha_{\langle\rangle}=\omega,\quad \alpha_{\langle i\rangle\conc\eta}=
\alpha^i_\eta.
  \end{array}\]
Checking that $\langle B^\ell_\eta,f^{\ell'}_\eta,\alpha_\eta:\eta\in T,\
\ell\leq\ell_\eta,\ \ell'\leq\ell_\eta'\rangle$ is as required is
straightforward. 

\noindent{\sc Case 3:}\qquad $\alpha=\beta+1$.\\
By (2) above and the proof of Case 2 we may assume that $B_m\in\ID^-_\alpha
(F)$. It follows from the definition of $\ID^-_\alpha(F)$ that there are
$f_\ell\in{}^{\textstyle\lambda}2$ (for $\ell\leq m$) such that
\[B^\oplus_\ell\stackrel{\rm def}{=}\{\delta<\lambda:\delta\mbox{ is limit
and }F(    f_ {\ell}   
\rest\delta)=0\Leftrightarrow\delta\in B_\ell\}\in\ID_\beta(F),\]
and hence $B^\oplus\stackrel{\rm def}{=}\bigcup\limits_{\ell\leq m}
B^\oplus_\ell\in\ID_\beta(F)$. Therefore $B^*_\ell\stackrel{\rm def}{=}
B_\ell\cap B^\oplus\in\ID_\beta(F)$. Now apply the inductive hypothesis for
$\beta$ and $\bar{B}^*=\langle B^*_\ell:\ell\leq m\rangle$ to get the
sequences $\langle B^{\ell,*}_\eta,f^{k,*}_\eta: \eta\in T^*,\ \ell\leq
\ell^*_\eta,\ k\leq k^*_\eta\rangle$ satisfying clauses (a)--(d) and such
that $\langle B^{\ell,*}_{\langle\rangle}:\ell\leq\ell^*_\eta\rangle=
\langle B^*_\ell:\ell\leq m\rangle$. Put 
\[\begin{array}{l}
T=\{\langle\rangle\}\cup\{\langle i\rangle:i<\lambda\}\cup\{\langle 0\rangle
\conc\eta:\eta\in T^*\},\\
\ell_{\langle 0\rangle\conc\eta}=\ell^*_\eta,\quad k_{\langle\rangle}=m+1,
\quad k_{\langle 0\rangle\conc\eta}=k_\eta,\\
B^\ell_{\langle 0\rangle\conc\eta}=B^{\ell,*}_\eta,\quad B^\ell_{\langle 0
\rangle\conc\langle i\rangle}=B_\ell\cap(i+\omega),\\
f^k_{\langle\rangle}=f_k,\quad f^k_{\langle
0\rangle\conc\eta}=f^{k,*}_\eta,\\ 
\alpha_{\langle\rangle}=\omega,\quad \alpha_{\langle 0\rangle\conc\eta}=
\alpha^*_\eta.
  \end{array}\]

\noindent (5)\quad If $f^\ell_\eta$ is not defined then choose $f^*$ as it.

\noindent 
(6), (7), (8), (9) Easy too.   
\end{proof}

\begin{remark}  \label{a11}
Note that it may happen that $\lambda\in\ida(F)$. However, if $\eta\in
{}^{\textstyle\lambda} 2$ is a weak diamond sequence for $F$ then the set
$\{\gamma<\lambda: \eta(\gamma)=0\}$ witnesses $\lambda\notin \ID_1^-(F)$. 
And conversely, if $\lambda\notin\ID_1^-(F)$ and $S^*\subseteq \lambda$
witnesses it, then the function $0_{S^*}\cup 1_{\lambda\setminus S^*}$ is a
weak diamond sequence for $F$.
\end{remark}

\begin{definition}  \label{a14}
\label{1.3Adef}
For a $\lambda$--colouring $F$ we define $\lambda$--colourings $F^\oplus$ and
$F^\otimes$ as follows.
\begin{enumerate}
\item A function $g\in{}^{\textstyle \gamma}(\Hla)$ is called
$F^\oplus$--standard if there is a tuple $(T,\bar{f},\bar{\alpha},\bar{A})$
(called a witness) such that 
\begin{enumerate}
\item[(i)] $T\subseteq {}^{\textstyle\omega{>}}\gamma$ is a well founded tree
(so $\langle\rangle\in T$, $\nu\vartriangleleft\eta\in T\ \Rightarrow \nu\in
T$ and $T$ has no $\omega$--branch);
\item[(ii)] $\bar{f}=\langle f_\eta^\ell:\eta\in T,\ \ell\leq k_\eta
\rangle$, where $f_\eta^\ell\in \DOM(F)\cap{}^{\textstyle\gamma}(\Hla)$;
\item[(iii)] $\bar{\alpha}=\langle\alpha_\eta:\eta\in T\rangle$, where
$\alpha_\eta<\lambda$;
\item[(iv)] $\bar{A}=\langle A^\ell_\eta:\eta\in T,\ \ell\leq\ell_\eta
\rangle$, where $A_\eta^\ell\subseteq\alpha_\eta$;
\item[(v)] $g(\beta)=(T\cap{}^{\textstyle\omega{>}}\beta,\langle f_\eta^\ell
\rest\beta:\eta\in T\cap{}^{\textstyle\omega{>}}\beta,\ \ell<k_\eta\rangle,
\langle\alpha_\eta:\eta\in T\cap{}^{\textstyle\omega{>}}\beta\rangle,\langle
A_\eta^\ell:\eta\in T\cap{}^{\textstyle\omega{>}}\beta,\ \ell\leq\ell_\eta
\rangle)$ for each $\beta<\gamma$.
\end{enumerate}
\item $\DOM(F^\oplus)=\bigcup\limits_{\alpha<\lambda}{}^{\textstyle\alpha}(
{\mathcal H}(\lambda))$ and for $g\in{}^{\textstyle\gamma}(\Hla)$:
\begin{enumerate}
\item[$(\oplus)_\alpha$] if $\gamma=0$ then $F^\oplus(g)=0$,
\item[$(\oplus)_\beta$]  if $\gamma>0$ and $g$ is not standard then
$F^\oplus(g)=0$, 
\item[$(\oplus)_\gamma$] if $\gamma>0$ and $g$ is standard as witnessed by
$\langle \bar{T},\bar{f},\bar{\alpha},\bar{A}\rangle$ then $F^\oplus(g)=
\bt^0_{F,g}(\langle\rangle)$, where $\bt^\ell_{F,g}(\eta)\in\{0,1\}$ (for
$\eta\in T$, $\ell=0,1$) are defined by downward induction as follows.
\begin{enumerate}
\item[If] $\eta\in\max(T)$ then $\bt^\ell_{F,g}(\eta)=1$ iff $\gamma\in
A_\eta$,
\item[if] $\eta\in T\setminus\max(T)$, $\gamma<\alpha_\eta$ then
$\bt^\ell_{F,g}(\eta)=1$ iff $\gamma\in A_\eta$,
\item[if] $\eta\in T\setminus \max(T)$, $\gamma\geq\alpha_\eta$ then
\[\begin{array}{lcl}
\bt^1_{F,g}(\eta)=1& \mbox{ iff }& F(f_\eta)=1\ \mbox{ or }\ (\exists i<
\gamma)(\bt^1_{F,g}(\eta\conc\langle i\rangle)=1),\\
\bt^0_{F,g}(\eta)=1& \mbox{ iff }& (\exists i<\gamma)(\bt^0_{F,g}(\eta\conc
\langle i\rangle)=1)\ \mbox{ or}\\
& & F(f_\eta   
)=1\ \&\ (\forall i<\gamma)(\bt^1_{F,g}(\eta\conc\langle i
\rangle)=0).
  \end{array}\]
\end{enumerate}
\end{enumerate}
\item A function $g\in {}^{\textstyle\gamma}(\Hla)$ is called
$F^\otimes$--standard if there is a tuple $(T,\bar{f},\bar{\ell},\bar{
\alpha},\bar{A})$ (called a witness) such that 
\begin{enumerate}
\item[(i)] $T\subseteq {}^{\textstyle\omega{>}}\gamma$ is a well founded tree;
\item[(ii)] $\bar{f}=\langle f_\eta:\eta\in T\rangle$, where $f_\eta\in\DOM(F)
\cap{}^{\textstyle\gamma}(\Hla)$;
\item[(iii)] $\bar{\ell}=\langle\ell_\eta:\eta\in T\rangle$, where
$\ell_\eta:{}^{\textstyle 3}\{0,1\}\longrightarrow\{0,1\}$;
\item[(iv)] $\bar{\alpha}=\langle\alpha_\eta:\eta\in T\rangle$, where
$\alpha_\eta<\lambda$;
\item[(v)] $\bar{A}=\langle A_\eta:\eta\in T\rangle$, where $A_\eta\subseteq
\alpha_\eta$;
\item[(vi)] $g(\beta)=(T\cap{}^{\textstyle\omega{>}}\beta,\langle f_\eta\rest
\beta:\eta\in T\cap{}^{\textstyle\omega{>}}\beta\rangle,\langle\ell_\eta:\eta
\in T\cap{}^{\textstyle\omega{>}}\beta\rangle,\langle\alpha_\eta:\eta\in T
\cap{}^{\textstyle\omega{>}}\beta\rangle,\langle A_\eta:\eta\in T
\cap{}^{\textstyle\omega{>}}\beta\rangle)$ for each $\beta<\gamma$.
\end{enumerate}
\item $\DOM(F^\otimes)=\bigcup\limits_{\alpha<\lambda}{}^{\textstyle\alpha}(
{\mathcal H}(\lambda))$ and for $g\in{}^{\textstyle\gamma}(\Hla)$:
\begin{enumerate}
\item[$(\otimes)_\alpha$] if $\gamma=0$ then $F^\otimes(g)=0$,
\item[$(\otimes)_\beta$]  if $\gamma>0$ and $g$ is not $F^\otimes$--standard
then $F^\otimes(g)=0$, 
\item[$(\otimes)_\gamma$] if $\gamma>0$ and $g$ is $F^\otimes$--standard as
witnessed by $\langle \bar{T},\bar{f},\bar{\ell},\bar{\alpha},\bar{A}\rangle$
then $F^\otimes(g)=\bt_{F,g}(\langle\rangle)$, where $\bt_{F,g}(\eta)\in\{0,
1\}$ (for $\eta\in T$) are defined by downward induction as follows. 
\begin{enumerate}
\item[If] $\eta\in\max(T)$ then $\bt_{F,g}(\eta)=1$ iff $\gamma\in A_\eta$,
\item[if] $\eta\in T\setminus\max(T)$, $1+\gamma<\alpha_\eta$ then
$\bt_{F,g}(\eta)=1$ iff $\gamma\in A_\eta$,
\item[if] $\eta\in T\setminus \max(T)$, $1+\gamma\geq\alpha_\eta$ then
\[\bt_{F,g}(\eta)=\ell_\eta(F(f_\eta),\max\{\bt_{F,g}(\eta\conc\langle\beta
\rangle):\beta<\gamma\},\min\{\bt_{F,g}(\eta\conc\langle\beta\rangle):
\beta<\gamma\}).\]
\end{enumerate}
\end{enumerate}
\end{enumerate}
\end{definition}

\begin{remark} 
On $ F^\oplus , F^\otimes $   
see \ref{a17}, \ref{a32}  below.
\end{remark} 

\begin{proposition}  \label{a17}
\label{1.3Bclaim}
Let $F$ be a $\lambda$--colouring. Then $F^\oplus$ is a $\lambda$--colouring
and 
\begin{enumerate}
\item[(a)] if $S\in\ida(F)$ then $0_S\cup 1_{\lambda\setminus S}\in\gB(
F^\oplus)$ and $\gB(F)\subseteq\gB(F^\oplus)$,
\item[(b)] $\ida(F)\subseteq\ID_1(F^\oplus)=\ID_1^-(F^\oplus)=\ida(F^\oplus)$, 
\end{enumerate}
\end{proposition}

\begin{proof}
(a)\quad Check.

\noindent (b)\quad  The main  
point is proving
 $\ida(F)\subseteq\ID_1(F^\oplus)$.\\

Suppose that $B\in\ida(F)$. We are going to show that then $B\in\ID^-_1(
F^\oplus)$. So suppose that $B'\subseteq B$. We want to find $g\in
\DOM_\lambda(F^\oplus)$ such that the set
\[\{\delta<\lambda:\delta\mbox{ is limit and } F(g\rest\delta)=0
\Leftrightarrow\delta\in B'\}\]
is in $\ID_0(F^\oplus)$ (what just means that it is non-stationary). Since
$B\in\ida(F)$ we have $B'\in\ida(F)$, so by \ref{easyprop}(4) we may find
$\langle B^\ell_\eta,f^k_\eta,\alpha_\eta:\eta\in T,\ \ell\leq\ell_\eta,\
k<k_\eta\rangle$ such that the clauses (a)--(d) of \ref{easyprop}(4) are
satisfied with $\ell_{\langle\rangle}=0$, $B'=B^0_{\langle\rangle}$. Define
$g$ as follows. For $\beta<\lambda$ let $T_\beta=T\cap{}^{\textstyle
\omega{>}}\beta$ and 
\[g(\beta)=\big(T_\beta,\ \langle f^k_\eta:\eta\in T_\beta, k\leq k_\eta
\rangle,\ \langle\alpha_\eta:\eta\in T_\beta\rangle,\ \langle B^\ell_\eta
\cap\alpha_\eta: \ell\leq\ell_\eta, \eta\in T_\beta\rangle\big).\]
Now look at the demands in \ref{1.3Adef}(2) -- they are exactly what
\ref{easyprop}(4) guarantees us.
\end{proof}

\begin{definition}   \label{a20}
\label{inequality}
Let $F_1,F_2$ be $\lambda$--colourings (with $\DOM(F_\ell)$ being either
${}^{\textstyle\lambda\!>}2$ or $\bigcup\limits_{\alpha<\lambda}
{}^{\textstyle\alpha}({\mathcal H}(\lambda))$, see \ref{colourings}(1)).
\begin{enumerate}
\item We say that $F_1\leq F_2$ \underline{when}  
there is $h: \DOM(F_1) \longrightarrow
\DOM(F_2)$ such that:   
\begin{enumerate}
\item[(a)] $\eta\trianglelefteq\nu\quad \Rightarrow\quad h(\eta)
\trianglelefteq h(\nu)$,
\item[(b)] $h(\eta)=\lim\limits_{\alpha<\delta} h(\eta\restriction\alpha)$,
for every $\eta\in {}^{\textstyle\delta}2$, $\delta$ a limit 
ordinal,   
\item[(c)] $(\forall\eta\in \DOM(F_1))(0<\lh(\eta)=\lh(h(\eta))\ \Rightarrow
\ F_1(\eta)=F_2(h(\eta)))$. 
\end{enumerate}
\item We say that $F_1\leq^* F_2$ \underline{when} 
there is $h:\DOM(F_1)\longrightarrow
\DOM(F_2)$ such that the clauses (a)--(c) above hold but
\begin{enumerate}
\item[(d)] if $\eta\in\DOM_\lambda(F_1)$ and $\lim\limits_{\alpha<\lambda}
h(\eta\rest\alpha)$ has length $<\lambda$ then $F_1(\eta\rest\alpha)=0$ for
every large enough $\alpha$.
\end{enumerate}
\end{enumerate}
\end{definition}

\begin{proposition} \label{a23}
\begin{enumerate}
\item $\leq^*$ and $\leq$ are transitive relations on $\lambda$--colourings,
satisfying   
$\leq^*\ \subseteq\ \leq$. 
\item $\leq$ is $\lambda^+$-directed.
\end{enumerate}
\end{proposition}

\begin{proposition}  \label{a26}
\label{prop1.8}
\begin{enumerate}
\item For every colouring $F_1:\bigcup\limits_{\alpha<\lambda}{}^{\textstyle
\alpha}({\mathcal H}(\lambda))\longrightarrow 2$ there is a colouring $F_2:
{}^{\textstyle\lambda\!>}2\longrightarrow 2$ such that $F_1\leq F_2\leq^*
F_1$.  
\item For every $\lambda$--colouring $F_2:{}^{\textstyle\lambda{>}}2
\longrightarrow 2$ there is a $\lambda$--colouring $F_1:\bigcup\limits_{
\alpha<\lambda}{}^{\textstyle\alpha}({\mathcal H}(\lambda))
\rightarrow 2   
$ such that
$F_2\leq F_1\leq^* F_2$.  
\end{enumerate}
\end{proposition}

\begin{proof}
1)\quad Let $F_1:\bigcup\limits_{\alpha<\lambda}{}^{\textstyle\alpha}({
\mathcal H}(\lambda))\longrightarrow 2$. Let $h_0$ be a one-to-one function
from $\Hla$ to ${}^{\textstyle\lambda{>}}2$, say $h_0(\eta)=\langle
\ell_{\eta,i}:i<\lh(h_0(\eta))\rangle$. Define a function
$h_1:\Hla\longrightarrow{}^{\textstyle\lambda{>}}2$ by: 
\[\begin{array}{l}
\lh(h_1(\eta))=2  
   \lh(h_0(\eta))+2,\\
h_1(\eta)(2i)=h_0(\eta)(i),\quad h_1(\eta)(2i+1)=0\quad\mbox{ for }i<\lh(
h_0(\eta)),\quad\mbox{ and}\\
h_1(\eta)(2\lh(h_0(\eta)))=h_1(\eta)(2\lh(h_0(\eta)+1))=1.
  \end{array}\]
Next, by induction on $\lh(\eta)$, we define a function $h^+:\bigcup\limits_{
\alpha<\lambda}{}^{\textstyle\alpha}(\Hla)\longrightarrow {}^{\textstyle
\lambda{>}}2$ as follows:
\[h^+(\langle\rangle)=\langle\rangle,\qquad h^+(\eta\conc\langle x\rangle)=
h^+(\eta)\conc h_1(x).\]

and if $ \eta \in \bigcup\limits_{\alpha<\lambda}{}^{\textstyle\alpha}({
\mathcal H}(\lambda)) $ has length the limit ordinal $ \delta $ 
then $ h^+( \eta ) = \cup  
 \{h^+(\eta \rest \beta ): \beta < \delta  \} $

Clearly $ h^+ $ is one to one with the right   
domain and range.

Finally we define a colouring $F_2:{}^{\textstyle\lambda{>}}2\longrightarrow
2$ by
\[F_2(\nu)=\left\{\begin{array}{ll}
F_1(\eta)&\mbox{if }\nu=h^+(\eta),\\
0        &\mbox{if }\nu\notin\rng(h^+).
		  \end{array}\right.\]

It is easy to check that $ F_2 $ is as required.  
\end{proof}

\begin{proposition}  \label{a29}
\label{anotprop}
Assume that $F_1,F_2$ are $\lambda$--colourings such that $F_1\leq F_2$, or  
$F_1\leq^* F_2$. \underline{Then}:   
\begin{enumerate}
\item For every $\eta\in    \DOM_ \lambda (F) $ 
there are $\nu\in \DOM_ \lambda (F) $ 
and a club $E$ of $\lambda$ such that 
\[(\forall\delta\in E)(F_1(\eta\rest\delta)=F_2(\nu\rest\delta)).\]
\item $\ID_\alpha(F_1)\subseteq\ID_\alpha(F_2)$, $\ID_\alpha^-(F_1)\subseteq
\ID_\alpha^-(F_2)$; hence $\ida(F_1)\subseteq\ida(F_2)$ and $\gB^+(F_1)
\subseteq\gB^+(F_2)$.
\item For every colouring $F$   
we have $\idc(F)\subseteq \wdmid  $ 
\end{enumerate}
\end{proposition}

\begin{proof}
Straightforward.
\end{proof}

\begin{conclusion}  \label{a32}
\label{con2}
Assume that $\lambda$ is a regular uncountable cardinal and $F:{}^{
\textstyle\lambda\!>}2\longrightarrow 2$ is a $\lambda$--colouring. Let
\[F^\otimes:\bigcup\limits_{\alpha<\lambda}{}^{\textstyle\alpha}({\mathcal
  H}(\lambda))\longrightarrow 2\]
be the colouring defined for $F$ in Definition \ref{1.3Adef}(4). 
Let $ \iota \in \{a,b\}$. 
\underline{Then}:
\begin{enumerate}
\item[(a)] $F\leq F^\otimes$.
\item[(b)] $\idi(F^\otimes)$ is a normal ideal on $\lambda$.
\item[(c)] $\gB(F)\subseteq\gB(F^\otimes)$ and $\idi(F)\subseteq\idi(
F^\otimes)=\wdmid(F^\otimes)$.
\item[(d)] $F^\otimes$ relates to itself as it relates to $F$, i.e.\quad 
{\em if}\quad $\alpha^*<\lambda^+$, $\langle S_\alpha:\alpha<\alpha^*\rangle$
is increasing continuous modulo $\idi(F^\otimes)$, $S_{\alpha+1}=S_\alpha\cup
A_\alpha\ \mod\ \idi(F^\otimes)$, $A_\alpha\in\gB(F^\otimes)$, $\ell_\alpha\in
2$, 

{\em then}\ \ for some $f\in {}^{\textstyle\lambda}({\mathcal H}(\lambda))$ 
\[\{\alpha<\lambda: F(f\rest\alpha)=1\}/{\mathcal D}_\lambda\]
is, in ${\mathcal P}(\lambda)/{\mathcal D}_\lambda$, the least upper bound
of the family $\{(A_\alpha\setminus S_\alpha)/{\mathcal D}_\lambda:
\ell_\alpha=1\}$ (where ${\mathcal D}_\lambda$ stands for the club filter). 
\item[(e)] The family $\gB(F^\otimes)$ is closed under complements, unions and
intersections of less than $\lambda$ sets, diagonal unions and diagonal
intersections and it includes bounded subsets of $\lambda$. Moreover
$\gB^+(F^\otimes)=\gB(F^\otimes)$. 
\item[(f)] If ${\mathcal P}(\lambda)/\idi(F^\otimes)$ is
$\lambda^+$--saturated \underline{then}  

for every set $X\subseteq\lambda$ there are sets $A,B\in 
   \mathfrak{B} 
(F^\otimes)$ such
that:
\begin{enumerate}
\item[$(\alpha)$] $A\subseteq X\subseteq B$,
\item[$(\beta)$]  for every $Y\in\gB(F^\otimes)$ one of the following occurs:
\begin{enumerate}
\item[(i)] \ \ \ the sets $(X\setminus A)\cap Y$, $(X\setminus A)\setminus
Y$, $(B\setminus X)\cap Y$, $(B\setminus X)\setminus Y$ are\footnote{
    hence
none of $X\setminus A$, $B\setminus A$ includes (modulo $\idi(F^\otimes)$) a
member of $\gB(F^\otimes)\setminus\idi(F^\otimes)$} not in $\idi(F^\otimes)$, 
\item[(ii)]\ \ $Y\cap (B\setminus A)\in\idi(F^\otimes)$,
\item[(iii)]\  $(B\setminus A)\setminus Y\in \idi(F^\otimes)$.
\end{enumerate}
\end{enumerate}
In the situation as above we denote $A=\max_{F^\otimes}(X)$, $B=\min_{
F^\otimes}(X)$ (note that these sets are unique 
only  
modulo $\idi(F^\otimes)$).
Moreover  
\item[(g)] In clause (f), 
if $A\subseteq\min_{F^\otimes}(B)$ then $\min_{F^\otimes}(A)
\subseteq\min_{F^\otimes}(B)\ \mod\ \idi(F^\otimes)$.
\item[(h)] In clause (f), when $ \iota = b $,
if $X\subseteq\lambda$, $X\notin\idi(F^\otimes)$ then for some
$Y_1,Y_2\subseteq X$ which are not in $\idi(F^\otimes)$ we have
\[{\max}_{F^\otimes}(Y_1)={\max}_{F^\otimes}(Y_2)=\emptyset\quad\mbox{ and
}\quad{\min}_{F^\otimes}(Y_1)={\min}_{F^\otimes}(Y_2)\notin\idi(F^\otimes).\]
    \item[(i)] 
    In clause (f), 
    $ \min _{F^ \otimes } $ and $ \max _{F^ \otimes }  $ commute 
    with the union of $ < \lambda $ and the intersection  of $< \lambda $ 
    sets.
\end{enumerate}
\end{conclusion}


\begin{proof}  
\underline{Clauses (a) and (b):}\qquad Should be clear.
\medskip

\underline{Clause (e):}\qquad Note that as $\theta=2$ we identify a sequence
$\eta\in {}^{\textstyle\lambda}2$ with $\{i<\lambda:\eta(i)=1\}$.
\smallskip

{\bf $\gB(F^\otimes)$ is closed under complementation.}\\
Suppose that $A\in\gB(F^\otimes)$. 
First, assume  
$A$ is bounded then let $g$,
$(T,\bar{f},\bar{\ell},\bar{\alpha},\bar{A})$ be as in \ref{1.3Adef}(3) with
$T=\{\langle\rangle\}\cup\{\langle i\rangle: i<\lambda\}$, $A_{\langle
\rangle}=\alpha_{\langle\rangle}\setminus A$, $\alpha_{\langle\rangle}>
\sup(A)$, $\ell_{\langle\rangle}$ constantly 1. Then $(\forall\alpha<
\lambda)(F^\otimes(g\rest(1+\alpha))=1\ \Leftrightarrow\ \alpha\in A)$, so
$F$ codes $\lambda\setminus A$. 

Second,  
suppose that $\sup(A)=\lambda$. Pick $g$
such that 
\[
(\forall\alpha<\lambda)(F^\otimes(g\rest (1+\alpha))=1\ \Leftrightarrow\
\alpha\in A).
\]
By our assumption, for arbitrarily large  $\beta<\lambda$ we have
$F^\otimes(g\rest\beta)=1$, so $g(\beta)$ is
\[
\big(T_\beta, \langle f^\beta_\eta:\eta\in T_\beta\rangle, \langle
\alpha^\beta_\eta:\eta\in T_\beta\rangle, \langle \ell^\beta_\eta:
\eta\in T_\beta\rangle, 
\langle A^\beta_\eta: \eta\in T_\beta\rangle\big)
\]
and it is as in \ref{1.3Adef}(3). If $\beta_1<\beta_2$ then the two values
necessarily cohere, in particular $T_{\beta_1}=T_{\beta_2}\cap{}^{\textstyle
\omega{>}}(\beta_1)$. Consequently there is $(T,\bar{f},\bar{\ell},\bar{
\alpha},\bar{A})$ such that $T=\bigcup\limits_{\beta<\lambda}T_\beta
\subseteq{}^{\textstyle\omega{>}}\lambda$ is closed under initial segments
and is well founded (as $T_\beta$ increase with $\beta$ and $\cf(\lambda)>
\aleph_0$). Thus we have proved
\begin{enumerate}
\item[$(\boxtimes)$] if $A\subseteq\lambda$ is unbounded and $F^\otimes$
coded by $g$ then there is ${\bf p}=(T,\bar{f},\bar{\ell},\bar{\alpha},
\bar{A})$ such that the clauses (i)--(vi) of \ref{1.3Adef}(3) hold for
$\gamma=\lambda$ and $g(\beta)={\bf p}\rest\beta$. 
\end{enumerate}
Now define ${\bf p}'$ like ${\bf p}$ (with the same $T$ etc) except that
$\ell^{{\bf p}'}_{\langle\rangle}=1-\ell^{{\bf p}}_{\langle\rangle}$ and
$A^{{\bf p}'}_{\langle\rangle}=A^{{\bf p}}_{\langle\rangle}$.
\smallskip

{\bf $\gB(F^\otimes)$ contains all bounded subsets of $\lambda$.}\\
By the first part of the arguments above all co-bounded subsets of $\lambda$
are in $\gB(F^\otimes)$, so (by the above) their complements are there too.
\smallskip

{\bf $\gB(F^\otimes)$ is closed under unions of length $<\lambda$.}\\
Let $B=\bigcup\limits_{i<\alpha}B_i$ where $\alpha<\lambda$ and $B_i\in
\gB(F^\otimes)$. Let $w=\{i<\alpha:\sup(B_i)=\lambda\}$ and for $i\in w$ let
$B_i$ be represented by $g_i\in{}^{\textstyle\lambda}(\Hla)$ which, by
$(\boxtimes)$, comes from ${\bf p}^i=(T^i,\bar{f}^i,\bar{\ell}^i,
\bar{\alpha}^i,\bar{A}^i)$. We may assume that $w=\beta\leq\alpha$. Let
\[\begin{array}{l}
T=\{\langle\rangle\}\cup\{\langle i\rangle: i<\lambda\}\cup\{\langle i
\rangle\conc\eta:\eta\in T^i,\ i<\beta\},\\
f_{\langle i\rangle\conc\eta}=f^i_\eta,\mbox{ etc}\\
\alpha_{\langle\rangle}\mbox{ is the first $\gamma\geq\omega$ such that }
\gamma\geq\alpha\ \&\ (\forall i\in [\beta,\alpha))(B_i\subseteq\gamma),\\
B_{\langle i\rangle}=\emptyset\quad\mbox{ if }i\geq\beta,\\
A_{\langle\rangle}=\bigcup\limits_{i<\alpha}B_i\cap\alpha_{\langle\rangle},\\
\ell_{\langle\rangle}(i_0,i_1,i_2)=i_1.
  \end{array}\]
Checking is straightforward.
\smallskip

{\bf $\gB(F^\otimes)$ is closed under diagonal unions.}\\
Let $B=\dun\limits_{i<\lambda}B_i$, where each $B_i\in\gB(F^\otimes)$ is
represented by $g_i\in{}^{\textstyle\lambda}(\Hla)$ which, by $(\boxtimes)$,
comes from ${\bf p}^i=(T^i,\bar{f}^i,\bar{\ell}^i,\bar{\alpha}^i,
\bar{A}^i)$. Let $T=\{\langle\rangle\}\cup\{\langle i\rangle\conc\eta:\eta
\in T_i,\ i<\lambda\}$, $f_{\langle i\rangle\conc\eta}=f^i_\eta$, etc,
$\alpha_{\langle\rangle}=\omega$, $B_{\langle\rangle}=B\cap\omega$ and
$\ell_{\langle\rangle}(i_0,i_1,i_2)=i_1$. 

So we have proved the first sentence in clause (e).  
The second sentence there follows by it and Def. \ref{a5}(2).  
Note that including the family of non-stationary sets     
follows by including the family of bounded subsets   
and being closed under diagonal unions.   
\medskip

\underline{Clause (c):}\qquad 
We concentrate on the case $ \idi= \ida$.
First note that $\gB(F)\subseteq\gB(F^\otimes)$ as
$\gB(F)\subseteq\gB^+(F)\subseteq\gB^+(F^\otimes)=\gB(F^\otimes)$ (the second
inclusion by (a) and \ref{anotprop}, the last equality by (e)). Next note
that 
\[\wdmid(F^\otimes)\subseteq\ID_1^-(F^\otimes)\subseteq\ID_1(F^\otimes)
\subseteq\idi(F^\otimes).\]
Now by induction on $\alpha$ we are proving that $\ID_\alpha(F^\otimes)
\subseteq \wdmid(F^\otimes)$. So suppose that we have arrived to a stage
$\alpha$.\\
If $\alpha=0$ then we use the fact that every non-stationary subset of
$\lambda$ is in $\gB(F^\otimes)$ (by (e)).\\
If $\alpha$ is limit then, by the induction hypothesis, $\ID^-_\alpha(
F^\otimes)\subseteq\gB(F^\otimes)$ and hence $\ID_\alpha
(F^  \otimes  )   
\subseteq\gB(
F^\otimes)$ (as $\gB(F^\otimes)$ is closed under diagonal unions by (e);
remember \ref{easyprop}(3)).\\
So suppose that $\alpha=\beta+1$ and $B\in\ID_\alpha(F^\otimes)$. Suppose
$B'\subseteq B$ (so $B'\in \ID_\alpha^-(F^\otimes)$). There is $B''\in\gB(
F)$ such that $B''\triangle B'\in\ID_\beta(F)$. By the first part we know
that $B''\in\gB(F^\otimes)$ and by the induction hypothesis $B'\triangle B''
\in\gB(F^\otimes)$. Consequently $B'\in\gB(F^\otimes)$.

Together we have proved that $\idi(F^\otimes)=\wdmid (F^\otimes)$. The
inclusion $\idi(F)\subseteq\idi(F^\otimes)$ is easy.

\underline{Clause (d)}
Easy.

\underline{Clause (f)}

Let ${  
{ {\mathscr D} }}_1 $ be  $ \{A \subseteq X: A \in 
\mathfrak{B} (F^ \otimes ) \setminus \idi (F^ \otimes )  
\} $ 
and let $ \{A_i: i < i_*  \} $ be a maximal sub-family of  
$ {{ {\mathscr D} }}_1 $   
such that 
$ i < j < i_* \Rightarrow A_i \cap A_j \in \idi(F^ \otimes ) $. 
By the assumption 
of clause (f) necessarily $ i_* < \lambda ^+$ so 
without lose of generality
$ i_* \le \lambda $. 
Let $  A $ be $ \cup \{A_i: i <  i_*  \} $ if $ i_* < \lambda $
and the diagonal union if $ i_* = \lambda $.
Clearly $ A \in \mathfrak{B} ^+(F)=\mathfrak{B} (F^ \otimes )  $. 

Let $ A' $ be chosen similarly replacing $ X $ by $ \lambda \setminus X$  and let
$ B = \lambda \setminus A'$

Clearly $ A , B $ are as required. 

\underline{Clause (g)}
Easy.


\underline{Clauses (h),  (i)}
 Easy.
\end{proof}

\begin{proposition}  \label{a35}
\label{1.3Cclaim}
Let $\lambda$ be a regular uncountable cardinal and $F$ be a
$\lambda$--colouring. 
\begin{enumerate}
\item If $\ID_\alpha(F)$ is $\lambda^+$--saturated then for some $\beta<
\lambda^+$ we have $\ID_{\alpha+\beta}(F)=\ida(F)$.
\item $\ID_\alpha(F)\subseteq\wdmid$,  
    see \ref{a5}(3); 
\item If $\ID_\alpha(F)$ is $\lambda^+$--saturated and $\lambda\notin\wdmid$
then $\wdmid=\ID_1(F')$ for some $\lambda$--colouring $F'$.
\item $\idb(F), \idc(F)$ are  normal ideals, and $\ID^1(F)\subseteq  \\
    \idb(F) \subseteq 
    \idc(F)\subseteq
\wdmid$.
\item $\ID^1(F^\otimes)=\wdmid(F^\otimes)$.
\item  $ \wdmid = \cup \{ \idi( F): F 
    $ a function form ${}^{ \lambda > } 2 $ to $  
    \}= 
     \cup \{\ID_1(F): F 
     $ a function form ${}^{ \lambda > } 2 $ to $
     \}  \cup \{\wdmid (F): F 
     $ a function form ${}^{ \lambda > } 2 $ to $
     \} $ 
     for $ \iota = a,b,c$.
\end{enumerate}
\end{proposition}

\begin{proof}
1)\qquad It follows from \ref{easyprop}(3) that $\ID_\gamma(F)$ increases
with $\gamma$
and $ \beta < \gamma , \ID_ \beta (F) = \ID_ {\beta + 1}$  
implyies   
$\ID_ \beta (F) = \ID_ \gamma  $; 
so the assertion should be clear.

2)\qquad By  the definition (and  \ref{con2}(c)).

3)\qquad Assume that $\ID_\alpha(F)$ is $\lambda^+$--saturated and
$\lambda\notin\wdmid$. By induction on $\beta<\lambda^+$ we try to 
    choose   
colourings $F_\beta$ such that
\begin{enumerate}
\item[(a)] $\ID(F) \subseteq \ID (F_  \alpha ) $ 
\item[(b)] if $\beta<\gamma$ then $\ida(F_\beta)\subseteq\ida(F_\gamma)$,
\item[(c)] $\ida(F_\beta)\neq\ida(F_{\beta+1})$.
\end{enumerate}
So we let $F_0=F$. If $\beta$ is limit then we use 
    \ref{a23}(2)   
to
choose $F_\beta$ so that $(\forall\gamma<\beta)(F_\gamma\leq F_\beta)$. 
Finally, if $\beta=\gamma+1$ then we let $F_\beta'=(F_\gamma)^\otimes$ (so
$\ida(F_\gamma)\subseteq\ID_1(F_\beta')=\ida(F_\beta')\subseteq\wdmid$). If
$\ida(F_\beta')\neq\wdmid$ then we choose a set $A\in\wdmid\setminus
\ida(F_\beta')$ and $F_\beta^*$ witnessing $A\in\wdmid$. We may assume that
$(\forall\alpha\in\lambda\setminus A)(\forall\eta\in{}^{\textstyle \alpha}2)
(F^*_\beta(\eta)=0)$. Now take a colouring $F_\beta$ such that $F_\beta',
F^*_\beta\leq F_\beta$.

After carrying out the construction choose $S^0_\beta\in\ida(F_{\beta+1})
\setminus\ida(F_\beta)$ (for $\beta<\lambda^+$) and let $S_\beta=S^0_\beta
\setminus\dun\limits_{\gamma<\beta}S^0_\gamma$. Then $\langle S_\beta:\beta<
\lambda^+\rangle$ is a sequence of pairwise disjoint members of ${\mathcal
P}(\lambda)\setminus\ida(F_0)\subseteq{\mathcal P}(\lambda)\setminus
\ID_\alpha(F)$, contradicting our assumptions. 

4), 5), 6) Easy too.   
\end{proof}

For the rest of this section we will assume the following

\begin{hypothesis}  \label{a38}
\label{tempassu}
\noindent (1)
We assume that 
\begin{enumerate}
\item[(a)] $\lambda$ is a regular uncountable cardinal, 
\item[(b)] $F$ is a $\lambda$--colouring, 
\item[(c)] $\lambda\notin\idb(F^\otimes)$, and 
\item[(d)] $\idb(F^\otimes)$ is $\lambda^+$--saturated, that is there is no
sequence $\langle A_\alpha:\alpha<\lambda^+\rangle$ such that for each
$\alpha<\beta<\lambda^+$ 
\[A_\alpha\notin\idb(F^\otimes)\quad\mbox{and}\quad\|A_\alpha\cap A_\beta\|<
\lambda.\]
\end{enumerate}

\noindent 
(2)
For each limit ordinal $\alpha\in [\lambda,\lambda^+)$ fix an enumeration
$\langle\varepsilon^\alpha_i:i<\lambda\rangle$ of $\alpha$.
\end{hypothesis}

\begin{construction}  \label{a41}
\label{1.4Aconstr}
Fix a sequence $\eta\in {}^{\textstyle\lambda} 2$ for a moment. We 
    choose  
a sequence 
\[\langle S_\alpha[\eta], A_\alpha[\eta], B_\alpha[\eta],\ell_\alpha[\eta],
m_\alpha[\eta], f_\alpha[\eta]:\alpha<\alpha^*[\eta]\rangle\]
as follows. By induction on $\alpha<\lambda^+$ we try to choose
$S_\alpha[\eta]=S_\alpha$, $A_\alpha[\eta]=A_\alpha$, $B_\alpha[\eta]=
B_\alpha$, $\ell_\alpha[\eta]=\ell_\alpha$, $m_\alpha[\eta]=m_\alpha$,
$f_\alpha[\eta]=f_\alpha$ such that: 
\begin{enumerate}
\item[(a)] $S_\alpha, A_\alpha, B_\alpha\subseteq\lambda$, $\ell_\alpha,
  m_\alpha\in \{0,1\}$, $f_\alpha\in {}^{\textstyle\lambda} 2$,
\item[(b)] $A_\alpha\notin\idb(F^\otimes)$, $A_\alpha\cap S_\alpha=\emptyset$,
\item[(c)] if $ \alpha =0 $ then $ S_\alpha = \emptyset $
\item[(d)]  
$S_{\alpha+1}=S_\alpha\cup A_\alpha$; 
\item[(e)] if $\alpha<\lambda$ is limit
  then $S_\alpha=\bigcup\limits_{\beta<\alpha} S_\beta $; if $\alpha\in
  [\lambda,\lambda^+)$ is limit then $S_\alpha=\{\gamma<\lambda: (\exists
  i<\gamma)(\gamma\in S_{\varepsilon^\alpha_i})\}$, $S_0=\emptyset$,
\item[(f)] $B_\alpha\in\idb(F^\otimes)$,
\item[(g)] for every $\delta\in\lambda\setminus(S_\alpha\cup B_\alpha)$
\[\eta(\delta)=m_\alpha\qquad\Rightarrow\qquad
F(f_\alpha\rest\delta)=\ell_\alpha,\]
\item[(h)] $A_\alpha=\{\delta\in\lambda\setminus S_\alpha: F(f_\alpha\rest
  \delta)=1-\ell_\alpha\}$. 
\end{enumerate}
It follows from \ref{tempassu} that at some stage $\alpha^*=\alpha^*[\eta]<
\lambda^+$ we get stuck (remember clause (b) above). Still, we may define 
$S_{\alpha^*}$ as in 
clause (c).
\end{construction}

\begin{proposition}  \label{a44}
\label{getF}
Assume \ref{tempassu}. Then:
\begin{enumerate}
\item There exists $\eta\in{}^{\textstyle\lambda}2$ such that 
\[\lambda\setminus S_{\alpha^*[\eta]}[\eta]\notin\idb(F^\otimes).\]
\item If $S\in\gB(F^\otimes)\setminus \idb(F^\otimes)$ then we can demand
$S\subseteq S_{\alpha^*[\eta]}[\eta]$.
\end{enumerate}
\end{proposition}

\begin{proof}
Assume not. Then for each $\eta\in{}^{\textstyle\lambda}2$ the set
$B_{\alpha^*}[\eta]\stackrel{\rm def}{=}\lambda\setminus S_{\alpha^*[\eta]}$
is in $\idb(F^\otimes)$. Now,
\[\{\alpha\in B_{\alpha^*}[\eta]:\eta(\alpha)=1\}\in\idb(F^\otimes)\subseteq
\gB(F^\otimes)\]
(see \ref{1.3Bclaim}). 

\begin{claim}  \label{a44b}
\label{clx}
For each $\alpha$, $S_\alpha\in\gB(F^\otimes)$. 
\end{claim}

\begin{proof}[Proof of the claim]
We show it by induction on $\alpha$. If $\alpha=0$ then $S_\alpha=\emptyset
\in\gB(F^\otimes)$ (see \ref{con2}(c)). If $\alpha<\lambda$ is a limit
ordinal then $S_\alpha=\bigcup\limits_{\beta<\alpha}S_\beta$ and by the
inductive hypothesis $S_\beta\in\gB(F^\otimes)$, so by \ref{con2}(e) we are
done (as $\gB(F^\otimes)$ is closed under unions of $<\lambda$ elements). If
$\alpha \in[\lambda,\lambda^+)$ is limit then we use the fact that $\gB(
F^\otimes)$ is closed under diagonal unions. If $\alpha=\beta+1$ then
$A_\beta\in\gB(F^\otimes )$ or $\lambda\setminus A_\beta\in\gB(F^\otimes )$ and hence we may
conclude that $A_\beta\in\gB(F^\otimes)$ (remember \ref{con2}(e)). Since
$\gB(F^\otimes)$ is closed under unions of length $<\lambda$ we are done.
\end{proof}

\begin{claim}  \label{a44d}
\label{cly}
For each $\alpha$, $Y_\alpha\stackrel{\rm def}{=}\{\beta<\lambda:\eta(\beta)
=1\}\cap S_\alpha\in\gB(F^\otimes)$.
\end{claim}

\begin{proof}[Proof of the claim]
We prove it by induction on $\alpha$. If $\alpha=0$ then
$Y_\alpha=\emptyset$ and there is nothing to do. The case of limit $\alpha$
is handled like that in the proof of \ref{clx}. So suppose that $\alpha=
\beta+1$. It suffices to show that the set $Y_\alpha\cap (S_\alpha\setminus
S_\beta)$ is in $\gB(F^ \otimes)$, which  
means that $Y_\alpha\cap A_\alpha$ is there
(remember clauses (g) and (h)
of \ref{a41}).    
Note that if $\delta\in A_\alpha\setminus
B_\alpha$ then $F(f_\alpha\rest\delta)=1-\ell_\alpha\neq\ell_\alpha$ and
hence $\eta(\delta)\neq m_\alpha$ so $\eta(\delta)=1-m_\alpha$. Consequently
$Y_\alpha\cap (A_\alpha\setminus B_\alpha)\in\{A_\alpha\setminus
B_\alpha,\emptyset\}$. But ${\mathcal P}(B_\alpha)\subseteq\gB(F^\otimes)$
so together we are done.
\end{proof}

It follows from \ref{clx}, \ref{cly} that 
\[\{\beta:\eta(\beta)=1\}\cap S_{\alpha^*[\eta]}[\eta]\in \gB(F^\otimes).\]
But $\lambda\setminus S_{\alpha^*[\eta]}[\eta]\in\idb(F^\otimes)$, so
${\mathcal P}(\lambda\setminus S_{\alpha^*[\eta]}[\eta])\subseteq\gB(
F^\otimes)$ so we get a contradiction.
\end{proof}

\begin{conclusion}  \label{a47}
\label{conc}
Assume \ref{tempassu}. Let $\eta\in{}^{\textstyle\lambda}2$, $X_\ell[\eta]=
(\lambda\setminus S_{\alpha^*[\eta]}[\eta])\cap\eta^{-1}(\{\ell\})$ (for
$\ell=0,1$). Then one of the following occurs:
\begin{enumerate}
\item[(A)] $\lambda\setminus S_{\alpha^*[\eta]}[\eta]\in\idb(F^\otimes)$,
\item[(B)] $X_0[\eta],X_1[\eta]\notin\idb(F^\otimes)$, and $X_0[\eta]\cup
X_1[\eta]\in \gB(F^\otimes)$, $X_0[\eta]\cap X_1[\eta]=\emptyset$, and for
every $f\in {}^{\textstyle\lambda}2$,
\end{enumerate}
{\em either}\quad the sequence $\langle F(f\rest\delta):\delta\in (\lambda
\setminus S_{\alpha^*[\eta]}[\eta])\rangle$ is $\idb(F^\otimes)$--almost
constant 

{\em or}\quad both sequences $\langle F(f\rest\delta):\delta\in X_0[\eta]
\rangle$ and $\langle F(f\rest\delta): \delta\in X_1[\eta]\rangle$ are not
$\idb(F^\otimes)$--almost constant.  
\end{conclusion}

\begin{proof}
Assume that the first possibility fails, so $\lambda\setminus
S_{\alpha^*[\eta]}[\eta]\notin\idb(F^\otimes)$.  

Assume $X_0[\eta]\in\idb(F^\otimes)$. Take any $f_{\alpha^*[\eta]}\in
{}^{\textstyle \lambda}2$ and choose $\ell_{\alpha^*[\eta]}\in\{0,1\}$ so
that
\[\{\delta\in\lambda\setminus S_{\alpha^*[\eta]}[\eta]: F(f_{\alpha^*[\eta]}
\rest\delta)=1-\ell_{\alpha^*[\eta]}\}\notin\idb(F^\otimes).\] 
Putting $m_{\alpha^*[\eta]}=0$ and $B_{\alpha^*[\eta]}=X_0[\eta]$ we get a
contradiction with the definition of $\alpha^*[\eta]$. Similarly one shows
that $X_1[\eta]\notin\idb(F^\otimes)$.

Suppose now that $f\in {}^{\textstyle\lambda}2$ and the sequence $\langle
F(f\rest\delta):\delta\in (\lambda\setminus S_{\alpha^*[\eta]}[\eta])
\rangle$ is not $\idb(F^\otimes)$--almost constant but, say, the sequence
$\langle F( f\rest\delta):\delta\in X_0[\eta]\rangle$ is
$\idb(F^\otimes)$--almost constant (and let the constant value be
$\ell_{\alpha^*[\eta]}$). Let $m_{\alpha^*[\eta]} =0$, 
$B_{\alpha^*[\eta]}=\{\delta\in X_0[\eta]: F(f\rest\delta)=1-\ell_{\alpha^*
[\eta]}\}$. Then $B_{\alpha^*[\eta]}\in\idb(F^\otimes)$ and since necessarily
\[\{\delta\in X_0[\eta]\cup X_1[\eta]: F(f\rest\delta)=1-\ell_{\alpha^*
[\eta]}\}\notin\idb(F^\otimes),\]
we immediately get a contradiction. Similarly in the symmetric case.
\end{proof}
  
\begin{remark}  \label{a49}
Note that,
in \ref{a47}, 
if $S\in\gB(F^\otimes)\setminus\idb(F^\otimes)$ then there is $\eta
\in{}^{\textstyle\lambda}2$ such that $\eta^{-1}[\{0\}]\supseteq\lambda
\setminus S$ and above $X_0,X_1\subseteq S$ and possibility (A) fails.
\end{remark}

\begin{proposition}    \label{a50}  
\label{1.7claim}
Assume \ref{tempassu}. 
\begin{enumerate}
\item We can find           
   $ S^*_F$, $S^*_0$ and $S^*_1$ such that:
\begin{enumerate}
\item[(a)] $S^*_F\in\gB(F^\otimes)$,
\item[(b)] $S^*_F=S^*_0\cup S^*_1$, $S^*_0\cap S^*_1=\emptyset$,  
    and $ S^*_0, S^*_1 $ witness $ S^*_F \in \idc_(F^\otimes ) $ 
\item[(c)] if $S^*_F\neq\lambda$ then $\idc(F^\otimes)\rest{\mathcal P}(
\lambda\setminus S^*_F)=\wdmid(F^\otimes)\rest{\mathcal P}(\lambda\setminus
S^*_F)$, $\lambda\setminus S^*_F\notin\idc(F^\otimes)$.
\item[(d)] if $S^*_F\neq\emptyset$ then $S^*_F\notin\idb(F^\otimes)$ and 
\[\{\big(S^*_0\cap F^\otimes(f)/\idb(F^\otimes),S^*_1\cap F^\otimes(f)/\idb(
F^\otimes)\big): f\in\DOM_\lambda\}\]
is an isomorphism from ${\mathcal P}(S^*_0)/\idb(F^\otimes)$ onto ${\mathcal
P}(S^*_1)/\idb(F^\otimes)$. 
\end{enumerate}
\item If in \ref{conc}, $S^*_F\subseteq     
S_{\alpha^*[\eta]}[\eta]\mod\idb(F)$
then we can add
\begin{enumerate}
\item[$(\circledast)$] for some $\rho\in{}^{\textstyle X_1}2$ for every $f\in
{}^{\textstyle\lambda}2$ we have
\[\{\delta\in X_1:F(f\rest\delta)=\rho(\delta)\}\neq\emptyset\mod\idb(
F^\otimes).\]
\end{enumerate}
\end{enumerate}
\end{proposition}

\begin{proof}
1)\qquad We try to choose by induction on $\alpha<\lambda^+$ sets
$S_\alpha,S_{\alpha,0},S_{\alpha,1}$ such that
\begin{enumerate}
\item[(a)] $S_\alpha\subseteq\lambda$,
\item[(b)] $S_\alpha=S_{\alpha,0}\cup S_{\alpha,1}$, $S_{\alpha,0}\cap
S_{\alpha,1}=\emptyset$, 
\item[(c)] if $\beta<\alpha$ and $\ell<2$ then
\[S_\beta\subseteq S_\alpha\mod\idb(F^\otimes)\quad\mbox{ and }\quad
S_{\beta,\ell}\subseteq S_{\alpha,\ell}\mod\idb(F^\otimes),\]
\item[(d)] the sets $S_{\alpha,0},S_{\alpha, 1}$ witness that $S_\alpha \in
\idc(F^\otimes)$ (see \ref{1.1Adef}(4)).
\end{enumerate}
At some stage $\alpha<\lambda^+$ we have to be stuck (as $\idb(F^\otimes)$ is
$\lambda^+$--saturated) and then $(S_\alpha,S_{\alpha,0},S_{\alpha,1})$ can
serve as $(S^*_F,S^*_0,S^*_1)$.

2)\qquad By the choice of $S^*_F$, for some $\ell<2$ we have
\[
{\mathcal P}(X_\ell)\neq\{F^\otimes(f)\cap X_\ell: f\in{}^{\textstyle
\lambda}2 \},   
\]
so let $Y\subseteq X_\ell$ be such that $Y\notin\{F^\otimes(f)\cap X_\ell:
f\in{}^{\textstyle \lambda} 2 \}$. Let $\rho= 0_Y\cup 1_{X_\ell\setminus
Y}$. Since without loss of generality $\ell=1$, we are done.
\end{proof}

\begin{remark}  \label{a53}
\begin{enumerate}
\item Recall that if $\lambda\notin\wdmid$ then $S^*_F\neq\lambda$.
\item  Recall:  $\ida(F^\otimes)   
)=\wdmid(F^\otimes)$ is a
normal ideal and $\idb (F^\otimes )\subseteq \idc(F^\otimes)$ are  
normal ideals   extending it.  
\end{enumerate}
\end{remark}



\section{Weak diamond for more colours}
In this section we deduce a weak diamond for, say, three colours, assuming 
the weak diamond for two colours and assuming that a certain ideal 
is saturated.

\begin{proposition}  \label{b2}
Assume that $\lambda$ is a regular uncountable cardinal and $\mu\leq 2^{<
\lambda}$. Let $F_i: {}^{\textstyle\lambda\!>}2\longrightarrow\{0,1\}$ be
$\lambda$--colourings for $i<\mu$. Then there is a colouring $F:{}^{
\textstyle\lambda>}2\longrightarrow\{0,1\}$ such that $F_i\leq F$ for every
$i<\mu$. 
\end{proposition}

\begin{proof}
{\sc Case 1.}\quad $\mu\leq 2^{\|\alpha\|}$ for some $\alpha<\lambda$.\\ 
Let $\rho_i\in {}^{\textstyle\alpha}2$ for $i<\mu$ be 
pairwise  
distinct. For $\eta\in
{}^{\textstyle\lambda\!>}2$ let $h_i(\eta)=\rho_i\conc\eta$. Define $F$ by:
\[F(\nu)=\left\{\begin{array}{ll}
0 &\mbox{if }\lh(\nu)<\alpha, \mbox{ or }\lh(\nu)\geq\alpha\\
\ &\mbox{but }\nu\restriction\alpha\notin\{\rho_i:i<\mu\},\\
F_i(\langle\nu(\alpha+\varepsilon):\varepsilon<\lh(\nu)-\alpha\rangle)& 
\mbox{if }\lh(\nu)\geq\alpha\mbox{ and for some  }i < \mu,\nu\restriction\alpha=\rho_i.
\end{array}\right. \]
It is easy to see that $F: {}^{\textstyle\lambda\!>}2\longrightarrow\{0,1\}$
and $h_i$ exemplifies that $F_i\leq F$.
\medskip

\noindent {\sc Case 2.}\quad $\mu=\lambda$.\\
For $\eta\in {}^{\textstyle\lambda\!>}2$, $i<\mu$, 
     we define $ h_i(\eta ) \in {}^{ i + 1 + \lh(\eta )}2$  as follows: for 
$\gamma <  i + 1 + \lh(\eta ) $ we 
let 
\[h_i(\eta)(\gamma)=\left\{\begin{array}{ll}
0 &\mbox{if }\gamma<i,\\
1 &\mbox{if }\gamma=i,\\
\eta(\gamma-(i+1))&\mbox{otherwise.}
\end{array}\right. \]
Next, for $\nu\in {}^{\textstyle\lambda\!>}2$ define:
\[F(\nu)=\left\{\begin{array}{ll}
F_i(\langle\nu(i+1+\gamma):\gamma<\lh(\nu)-(i+1)\rangle)&\mbox{if }i=\min\{j:
\nu(j)=1\}\\
0&\mbox{if there is no such $i$.} 
\end{array}\right. \]
Now check.
\medskip

\noindent {\sc Case 3.}\quad Otherwise, for each $\alpha<\lambda$ choose
$F^{\alpha}:{}^{\textstyle\lambda\!>}2\longrightarrow\{0,1\}$ such that
$(\forall i<2^{\|\alpha\|})(F_i\leq F^\alpha)$ (exists by Case 1). Let 
$F:{}^{\textstyle\lambda\!>}2\longrightarrow\{0,1\}$ be such that $(\forall
\alpha<\lambda)(F^\alpha\leq F)$ (exists by Case 2). 

The proposition follows. 
\end{proof}

\begin{theorem}   \label{b5}   
\label{2.6thm}
If (A) then (B) where  
    \begin{enumerate} 
\item[(A)]  
\begin{enumerate} 
    \item[(a)] $\lambda$ is a regular uncountable cardinal.  
    \item[(b)] $F^{\tr}:
{}^{\textstyle\lambda\!>}2\longrightarrow 3$. 
    \item[(c)] For $i<3$ let $F_i:
{}^{\textstyle\lambda\!>}2\longrightarrow \{0,1\}$ be such that  
\[F_i(\eta)=\left\{\begin{array}{ll}
1&\mbox{if } F^{\tr}(\eta)=i,\\
0&\mbox{otherwise, }
\end{array} \right. \],
    \item[(d)] 
    Let $F: {}^{\textstyle\lambda\!>}2\longrightarrow\{0,1\}$ be such that
$(\forall i<3)(F_i\leq F)$. Moreover 
$ F'_2 \le F$ where $ F'_ 2 $ is the function with domain $ {}^{ \lambda < } 2 $ defined by
$ F'_2(\eta ) = \min  \{1, F( \eta ) \} $.
(remember \ref{anotprop}(3)), 
    \item[(e)] $\idb(F^\otimes)$ is
$\lambda^+$--saturated, i.e.\ there is no sequence
 $\langle A_\alpha:\alpha <\lambda^+\rangle$ such 
 that\footnote{
    As is well known, writing below 
        $ A_ \alpha \cap A_ \beta \in \ID^b(F^\otimes )$ 
        instead $ \| A_\alpha \cap A_ \beta \|  < \lambda $
    does not change anything.
    }  
\[
(\forall\alpha<\beta<\lambda^+)(A_\alpha\notin\idb(F^\otimes )\quad\&\quad\|
A_\alpha\cap A_\beta\|<\lambda).
\]
\end{enumerate} 
\item[(B)] 
\underline{Then} there is a weak diamond sequence for $F^{\tr}$, even for every $S\in
\gB(F^\otimes)\setminus\idc(F^\otimes)$.
\end{enumerate} 
\end{theorem}

\begin{proof}
Let $S^*_F$ be as in \ref{1.7claim}. Since 
    we are assuming 
$\lambda\notin\idc(F^\otimes)$
necessarily $\lambda\setminus S^*_F\notin\idb(F^\otimes)$.   

[Why? because   by \ref{a53}(1)(b) we have 
$\idc(F^\otimes)=\idb(F^ \otimes )+S^*_F$].   

It follows from \ref{getF} and \ref{conc} that there are disjoint sets $X_0,
X_1\subseteq\lambda$ (even disjoint from $S^*_F$ from \ref{1.7claim}) such
that $X_0,X_1\notin\idb(F^\otimes)$, $X_0\cup X_1\in\gB(F^\otimes)$ and for
every $f\in {}^{\textstyle\lambda}2$ we have one of the following: 
\begin{enumerate}
\item[(a)] the sequence $\langle F(f\rest\delta):\delta\in X_0\cup X_1\rangle$
is $\idb(F^\otimes)$--almost constant, or
\item[(b)] both sequences $\langle F(f\rest\delta):\delta\in X_0\rangle$ and
$\langle F(f\rest\delta): \delta\in X_1\rangle$ are not $\idb(F^\otimes
)$--almost constant.  
\end{enumerate}
It follows from \ref{1.7claim}(2) that we may assume that there is $\eta\in
{}^{\textstyle X_1}2$ such that for every $f\in {}^{\textstyle \lambda}2$
the set
\[\{\delta\in X_1:F(f\rest\delta)=\eta(\delta)\}\]
is stationary. Define a function $\rho \in {}^{\textstyle\lambda }2$ as
follows:  
\[\rho(\alpha)=\left\{
\begin{array}{ll}
1+\eta(\alpha)&\mbox{if } \alpha \in X_1,\\
0&\mbox{otherwise. }
\end{array} \right. \]

\begin{claim}   \label{b5c}
\label{cl2}
$\rho$ is a weak diamond sequence for $F^\tr$ even on $X_0\cup X_1$.
\end{claim}

\begin{proof}[Proof of the claim]
Let $f\in {}^{\textstyle\lambda}2$
    and we shall prove that $ Y_f = \{\delta \in X_0 \cup X_1  
    :    
    F^{\tr}(\eta \rest \delta )= \eta (\delta )\} 
    \not=\emptyset  \mod $  $ \idb (F ^ \otimes ) $.   

If $\{\alpha\in X_0: F^\tr(f\rest
\alpha)=2 \}\notin\idb(F)$ then we are done (remember \ref{easyprop}(3)). 
Otherwise 
    (by the definition of $ F_0 $), 
we have 
\[\{\alpha \in X_0: F_2 (f\restriction\alpha)=1\}\in \idb(F).\]
For $\ell<3$,  
    as $ F_ {\ell} \le F$  
let $f_\ell\in{}^{\textstyle\lambda}2$ be such that the set
$\{\alpha<\lambda: F_\ell(f\rest\alpha)=F(f_\ell\rest\alpha)\}$ contains a
club of $\lambda$ 
and $ g \in {}^{ \lambda } 2 $  such that the set
$ \{\alpha < \lambda : F'_2(g \rest \alpha )= F(f \rest \alpha ) \} $ 
contains a club of $ \lambda $, exist by clause (A)(d) of the assumption of the theorem.
 
We  now     
  use $f_2$.  
Then
\[\{\alpha\in X_0: F(f_2\rest\alpha)=1\}\in\idb(F^\otimes),\]
and hence, by the choice of the sets $X_0,X_1$,  
    clause (b) there fails hence clause (a) holds, so 
\[\{\alpha\in  X_1 \cup   
    X_1: F(f_2\rest\alpha)=1\}\in\idb(F^\otimes).\]
Consequently, 
\[Z = 
   \{\alpha\in X_1: F^\tr(f\rest\alpha)=2   
   \}=\{\alpha\in X_1: F_2
   (f\rest  \alpha)=1\}\in\idb (F^\otimes).\]
   
Now we use the choice of $\eta$,
by it we 
know that the set
\[ Y=\{\delta\in X_1: F(g  
\rest\delta)=\eta(\delta)\} \]
is stationary and even $ \not=  \emptyset \mbox{ mod } \idb(F^ \otimes )$.  
Hence for some $k\in \{0,1\}$ the set
\[Y_k=\{\delta\in X_1: F(g   
\rest\delta)=k=\eta(\delta)\}\]
is stationary and even $ \not=  \emptyset \mbox{ mod }\idb(F^ \otimes )$, but 
$\{\delta\in X_1: F(g  
\rest\delta)=F'_2(f\rest\delta)\}$
contains a club. Hence 
\[Y^*_k=\{\delta\in X_1: F(g   
\rest\delta)=k=\eta(\delta)\ \mbox{ and }\
F(g  
    \rest\delta)=F'_2(f\rest\delta) \mbox{ and }  \delta \notin Z  
      \}\] 
is stationary and even $ \not=  \emptyset  \mbox{ mod }   \idb(F^ \otimes ) $.
    Finally note that    
\[\delta\in Y^*_k\quad\Rightarrow\quad F(f_1\rest\delta)=\eta(\delta)=
F_1(f\rest\delta)=k\quad\Rightarrow\quad F^{\tr}(f\rest\delta)=k.\]
Thus the   
claim and the theorem are proved.
\end{proof}
\end{proof} 

\begin{theorem}   \label{b8}
\label{claim3}
Suppose $F^\tr$ is a $(\lambda,\theta)$--colouring, $\theta\leq\lambda$ and
$F_i$ (for $i<\theta$) are given by
\[F_i(f)=\left\{\begin{array}{ll}
1&\mbox{if }F(f)=i,\\
0&\mbox{otherwise.}
\end{array}\right.\]
Let $F:{}^{\textstyle\lambda\!>}2\longrightarrow 2$ be such that $(\forall
i<\theta)(F_i\leq F)$ and let $F^\otimes$ be as in \ref{1.3Adef} for $F$. 
Suppose that $\idb(F^\otimes)$ is $\lambda^+$--saturated, and $S^*_{F^\otimes}
\neq\lambda$  
(equivalently $\lambda\notin\idc(F^\otimes)$).

If $ (\otimes )$  then $( \bigstar)  $  \underline{where}:

\begin{enumerate}
\item[$(\otimes)$] there are sets $Y_i\subseteq\lambda\setminus S^*_{
F^\otimes}$ for $i<\theta$ such that 
\begin{enumerate}
\item[(a)] $(\forall i<\theta)(Y_i\notin\idb(F^\otimes))$,
\item[(b)] the sets $Y_i$ are pairwise disjoint or at least 
\[(\forall i<j<\theta)(Y_i\cap Y_j\in\idb(F^\otimes)),\]
\item[(c)] $\bigcap\limits_{i<\theta}\min_{F^\otimes}(Y_i)\notin
\idb(F^\otimes)$, see \ref{con2}(h). 
\end{enumerate}
\end{enumerate}

\begin{enumerate}
\item[$(\bigstar)$] there is a weak diamond sequence $\eta\in {}^{\textstyle
\lambda}\theta$ for $F^\tr$, i.e.
\[
(\forall f\in {}^{\textstyle\lambda}2)(\{\delta<\lambda: F^\tr(f\rest
\delta)=\eta(\delta)\} 
   \mbox{ is stationary    and even }   \not=  \emptyset 
  \mod  
  \idb(F^ \otimes ) );
\]
moreover
\[(\forall f\in {}^{\textstyle\lambda}2)(\{\delta<\lambda: F^\tr(f\rest
\delta)=\eta(\delta)\}\notin\idb(F^\otimes)).\]
\end{enumerate}
\end{theorem}

\begin{proof}
We may assume that the sets $\langle Y_i:i<\theta\rangle$ are pairwise
disjoint (otherwise we use $Y'_i=Y_i\setminus\bigcup\limits_{j<i} Y_j$). Let
$\eta\in{}^{\textstyle\lambda}\theta$ be such that $(\forall i<\theta)(\eta
\rest Y_i=i)$. Note that if
\[\{\delta\in Y_i: F^\tr(f\rest \delta)=i\}\in\idb(F^\otimes)\]
then we also have
\[\{\delta<\lambda: F^\tr(f\rest\delta)=i\}\in\gB(F^\otimes)\]
(use $F_i\leq F\leq F^\otimes$). Consequently, in this case, we have
\[\{\delta\in{\min}_{F^\otimes}(Y_i): F^{\tr}(f\rest\delta)=i\}\in\idb(
F^\otimes).\]
If this occurs for every $i<\theta$ then
\[\{\delta\in\bigcap_{i<\theta}{\min}_{F^\otimes}(Y_i): (\exists
i<\theta)(F(f\rest\delta)=i)\}\in\idb(F^\otimes),\]
but for each $\delta$, for some $i<\theta$ we have $F(f\rest\delta)=i$, a
contradiction. 
\end{proof}

\begin{proposition}   \label{b11}
\label{claim5}
Under the assumptions of 
    \ref{b8}   
(so the ideal $\idb(F^\otimes)$ is
$\lambda^+$--saturated), if $X\subseteq\lambda\setminus S^*_{F^\otimes}$,
$X\notin\idb(F^\otimes)$ then there is a partition $(X_0,X_1)$ of $X$ (so
$X_0\cup X_1=X$, $X_0\cap X_1=\emptyset$) such that 
\[X_0,X_1\notin\idb(F^\otimes),\quad\mbox{and}\quad {\min}_{F^\otimes}(X_0)=
{\min}_{F^\otimes}(X_1)={\min}_{F^\otimes}(X).\] 
\end{proposition}

\begin{proof}
Let
\[\begin{array}{ll}
{\mathcal A}_{F^\otimes}\stackrel{\rm def}{=}\{Z\subseteq\lambda:&
Z\notin\idb(F^\otimes)\mbox{ and there is a partition $(Z_0,Z_1)$ of $Z$}\\
\ &\mbox{such that }{\min}_{F^\otimes}(Z_0)={\min}_{F^\otimes}(Z_1)\ \mod\
\idb(F^\otimes)\}. 
\end{array}\]
Note that, by \ref{con2}(h),
\begin{enumerate}
\item[$(*)$] \qquad $(\forall Y\in\idb(F^\otimes)^+)(\exists Z\in {\mathcal
A}_{F^\otimes})(Z\subseteq Y)$.
\end{enumerate}
Let $X\subseteq\lambda$, $X\notin\idb(F^\otimes)$ and let $\langle Z_\alpha:
\alpha<\alpha^*\rangle$ be a maximal sequence such that for each $\alpha<
\alpha^*$:
\[Z_\alpha\in {\mathcal A}_{F^\otimes},\quad Z_\alpha\subseteq X,\quad\mbox{
and}\quad(\forall\beta<\alpha)(Z_\alpha\cap Z_\beta\in\idb(F^\otimes)).\] 
Necessarily $\alpha^*<\lambda^+$, so without loss of generality $\alpha^*\leq
\lambda$, $\min(Z_\alpha)>\alpha$ and $Z_\alpha\cap Z_\beta=\emptyset$ for
$\alpha<\beta<\alpha^*$. Let $\langle Z^0_\alpha, Z^1_\alpha\rangle$ be a
partition of $Z_\alpha$ witnessing $Z_\alpha\in {\mathcal A}_{F^\otimes}$. Put
\[Z_0\stackrel{\rm def}{=}\bigcup_{\alpha<\alpha^*} Z^0_\alpha\quad\mbox{ and
}\quad Z_1\stackrel{\rm def}{=}\bigcup_{\alpha<\alpha^*} Z^1_\alpha.\]
Then $Z_0\cap Z_1=\emptyset$, $Z_0\cup Z_1\subseteq X$. Note that
$\bigcup\limits_{\alpha<\alpha^*} Z_\alpha$ is equal to the diagonal union
and, by $(*)$ above, $X\setminus\bigcup\limits_{\alpha<\alpha^*} Z_\alpha\in
\idb(F^\otimes)$. Consequently we may assume $Z_0\cup Z_1= \bigcup\limits_{
\alpha<\alpha^*} Z_\alpha=X$. Next, since 
\[{\min}_{F^\otimes}(Z_0)\supseteq{\min}_{F^\otimes}(Z^0_\alpha)\supseteq
Z^0_\alpha\cup Z^1_\alpha=Z_\alpha,\]
we get
\[{\min}_{F^\otimes}(Z_0)\supseteq\bigcup_{\alpha<\alpha^*} Z_\alpha=X=Z_0\cup
Z_1,\]
and similarly one shows that $\min_{F^\otimes}(Z_1)\supseteq X$. Now we use
\ref{con2}(h) to finish the proof. 
\end{proof}

\begin{proposition}   \label{b14}
Under the assumptions of \ref{claim3}:  
\begin{enumerate}
\item If $2^\theta<\lambda$ then there is a sequence $\langle Y_i:
i<\theta\rangle$ as required in \ref{claim3}$(\otimes )$
\item Similarly if $\theta\leq\aleph_0$.
\item In both cases, if $S\notin\idb(F^\otimes)$ then we can demand $(\forall
i<\theta)(Y_i\subseteq S)$.
\end{enumerate}
\end{proposition}

\begin{proof}
1)\quad By induction on $\alpha\leq\theta$ we choose sets $X_\eta
\subseteq \lambda$ for $\eta\in {}^{\textstyle\alpha}2$ such that:
\begin{enumerate}
\item[(i)]\ \ \ $X_{\langle\rangle}\notin\idb(F^\otimes)$,\
\item[(ii)]\ \  if $\alpha$ is limit then $X_\eta=\bigcap\limits_{i<\alpha}
X_{\eta\rest i}$,
\item[(iii)] if $\alpha=\beta+1$, $\eta\in {}^{\textstyle\beta}2$ and
$X_\eta\in \idb(F^\otimes)$ then $X_{\eta\conc\langle 0\rangle}=X_\eta$,
$X_{\eta\conc\langle 1\rangle}=\emptyset$;

if $\alpha=\beta+1$, $\eta\in {}^{\textstyle\beta}2$ and
$X_\eta\notin \idb(F^\otimes)$ then $(X_{\eta\conc\langle 0\rangle},X_{\eta
\conc\langle 1\rangle})$ is a partition of $X_\eta$ such that $\min_{
F^\otimes}(X_{\eta\conc\langle 0\rangle})=\min_{F^\otimes}(X_{\eta\conc
\langle 1\rangle})=\min_{F^\otimes}(X_\eta)$. 
\end{enumerate}
It follows from \ref{claim5} that we can carry out the construction.

Clearly $\langle X_\eta:\eta\in{}^{\textstyle\theta}2\rangle$ is a partition
of $X_{\langle\rangle}$, so (as $2^\theta<\lambda$ and $\idb(F^\otimes)$ is
$\lambda$--complete) we can find a sequence $\eta\in{}^{\textstyle\theta}2$
such that $X_\eta\notin\idb(F^\otimes)$. Then 
\[(\forall\alpha<\theta)(X_{\eta\rest\alpha}\notin\idb(F^\otimes))\]
(as each of these sets includes $X_\eta$). Moreover, for each $\alpha<\theta$
and for $\ell=0,1$ we have 
\[{\min}_{F^\otimes}(X_{\eta\rest\alpha\conc\langle\ell\rangle})\supseteq
X_{\eta\rest \alpha}\supseteq X_\eta.\]
Put $Y_\alpha=X_{\eta\rest\alpha\conc\langle 1-\eta(\alpha)\rangle}$. Then
$\langle Y_\alpha:\alpha<\theta\rangle$ is a sequence of pairwise disjoint
sets (as $X_{\eta\rest\alpha\conc\langle 0\rangle}\cap X_{\eta\rest\alpha
\conc\langle 1\rangle}=\emptyset$) and for every $\alpha<\theta$
\[Y_\alpha\notin\idb(F^\otimes)\quad\mbox{ and }\quad {\min}_{F^\otimes}(
Y_\alpha)\supseteq X_{\eta\rest\alpha}\supseteq X_\eta.\]
Hence $\bigcap\limits_{\alpha<\theta}\min_{F^\otimes}(Y_\alpha)\notin\idb(
F^\otimes)$. Let $Z_\alpha=Y_\alpha\cap\min_{F^\otimes}(X_\eta)$. Note that
$\min_{F^\otimes}(Z_\alpha)=\min_{F^\otimes}(X_\eta)$ (the ``$\leq$'' is
clear; if $\min_{F^\otimes}(Z_\alpha)<\min_{F^\otimes}(X_\eta)$ then
$\min_{F^\otimes}(X_\eta)\setminus\min_{F^\otimes}(Z_\alpha)$ contradicts the
definition of $\min_{F^\otimes}(Y_\alpha)$). Thus the sequence $\langle
Z_\alpha:\alpha<\theta\rangle$ is as required. Moreover 
\[{\min}_{F^\otimes}(Z_\alpha)=\bigcup\limits_{\beta}{\min}_{F^\otimes}(
Z_\beta).\]

\noindent 2)\quad Let $X\subseteq\lambda$, $X\notin\idb(F^\otimes)$. By
induction on $n$ we choose sets $X_n',X''_n$ such that $X'_n\cap X''_n=
\emptyset$, $X'_n\cup X''_n\supseteq X$, and
\[{\min}_{F^\otimes}(X'_n)={\min}_{F^\otimes}(X''_n)={\min}_{F^\otimes}(X).\]
For $n=0$ we use \ref{claim5} for $X$ to get $X'_0,X''_0$. For $n+1$ we use
\ref{claim5} for $X''_n$ to get $X'_{n+1},X''_{n+1}$.

Finally we let $Y_n=X''_n$ (note that 
$\min_{F^\otimes}(Y_n)=\min_{F^\otimes}(X)$). 
\end{proof}

\begin{conclusion}   \label{b17}  
Assume that
\begin{enumerate}
\item[(A)] $\lambda$ is a regular uncountable cardinal,
\item[(B)] $F$ is a $(\lambda,\theta)$--colouring such that $\lambda\notin
\idb(F)$ and $\idb(F)$ is $\lambda^+$--saturated,
\item[(C)] $2^\theta<\lambda$ or $\theta=\aleph_0$,
\item[(D)] $(\exists\mu<\lambda)(2^\mu=2^{{<}\lambda}<2^\lambda)$ or at least
$\lambda\notin\wdmid$ or at least $\lambda\notin\idc(F)$.
\end{enumerate}
Then there is a weak diamond sequence for $F$. Moreover, there is $\eta\in
{}^{\textstyle \lambda}\theta$ such that for each $f\in\DOM_\lambda(F)$ we
have
\[\{\delta<\lambda: F(f\rest\delta)=\eta(\delta)\}\notin\idb(F).\]
\end{conclusion}

\section{An application of Weak Diamond}
In this section we present an application of Weak Diamond in model
theory. For more on model--theoretic investigations of this kind we refer
the reader to \cite{Sh:576} and earlier work \cite{Sh:88}
 (and see \cite{Sh:88r}), and to an
excellent survey my Makowsky, \cite{Mw85}.
 
\begin{definition}   \label{c2}  
\label{amalgdef}
Let $\gK$ be a collection of models.
\begin{enumerate}
\item For a cardinal $\lambda$, $\gK_\lambda$ stands for the collection of all
members of $\gK$ of size $\lambda$.
\item We say that a partial order $\leq_\gK$ on $\gK_\lambda$ is
$\lambda$--nice if 
\begin{enumerate}
\item[($\alpha$)] $\leq_\gK$ is a suborder of $\subseteq$ and it is closed
under isomorphisms of models (i.e.\ if $M,N\in\gK_\lambda$, $M\leq_\gK N$
and $f:N\longrightarrow N'\in\gK_\lambda$ is an isomorphism then $f[M]
\leq_\gK N'$), 
\item[($\beta$)]  $(\gK_\lambda,\leq_\gK)$ is $\lambda$--closed (i.e.\ any
$\leq_\gK$--increasing sequence of length $\leq\lambda$ of elements of
$\gK_\lambda$ has a $\leq_\gK$--upper bound in $\gK_\lambda$) and
\item[($\gamma$)] if $\bar{M}=\langle M_\alpha: \alpha<\lambda\rangle$ is an
$\leq_\gK$--increasing sequence of elements of $\gK_\lambda$ then
$\bigcup\limits_{\alpha<\lambda} M_\alpha$ is the $\leq_\gK$--upper bound to
$\bar{M}$ (so $\bigcup\limits_{\alpha<\lambda} M_\alpha\in\gK_\lambda$). 
\end{enumerate}
\item Let $N\in\gK_\lambda$, $A\subseteq |N|$. We say that {\em the pair
$(A,N)$ has the w.  amalgmation property in $\gK_\lambda$} if for every $N_1,N_2
\in\gK_\lambda$ such that $N\leq_\gK N_1$, $N\leq_\gK N_2$ there are $N^*\in
\gK_\lambda$ and $\leq_{\gK}$--embeddings $f_1,f_2$ of $N_1,N_2$ into $N^*$,
respectively, such that $f_1( A)= f_2( A)$. (In words: $N_1,N_2$ can
be amalgamated over $(A,N)$ setwise.)
\item We say that $(\gK,\leq_\gK)$ has {\em the amalgamation property for
$\lambda$} if for every $M_0,M_1,M_2\in\gK_\lambda$ such that $M_0\leq_{\gK}
M_1$, $M_0\leq_{\gK} M_2$ there are $M\in\gK_\lambda$ and
$\leq_{\gK}$--embeddings $f_1,f_2$ of $M_1,M_2$ into $M$, respectively, such
that 
\[M_0\leq_{\gK} M\quad\mbox{ and }\quad f_1\rest M_0=f_2\rest M_0=\id_{M_0}.\]
\end{enumerate}
\end{definition}

\begin{theorem}   \label{c5}
Assume that $\lambda$ is a regular uncountable cardinal for which the weak
diamond holds (i.e.\ $\lambda\notin\wdmid$). Suppose that $\gK$ is a class of
models, $\gK$ is categorical in $\lambda$ (i.e.\ all models from $\gK_\lambda$
are isomorphic), it is closed under isomorphisms of models, and $\leq_\gK$ is
a $\lambda$--nice partial order on $\gK_\lambda$ and $M\in\gK_\lambda$. Let
$\bar{A}=\langle A_\alpha: \alpha<\lambda\rangle$ be an increasing continuous
sequence of subsets of $|M|$ such that
\[(\forall\alpha<\lambda)(\|A_\alpha\|<\lambda)\quad\mbox{ and }\quad
\bigcup_{\alpha<\lambda} A_\alpha=M.\]
Then the set
\[S^{\bar{A}}_M\stackrel{\rm def}{=}\big\{\alpha<\lambda\!: (A_\alpha,M)\mbox{
does not have the w.  amalgmation property}\/\big\}\]
is in $\wdmid$.
\end{theorem}
 
\begin{proof} Assume that $S^{\bar{A}}_M\notin\wdmid$. 

We can fix a partition $ \langle D_i: i < \lambda \rangle $  of $ \lambda $ 
    to sets each of cardinality $ \lambda $.
    
We may assume that $|M|=\lambda$. By induction on $i<\lambda$ we choose pairs
$(B_\eta,N_\eta)$ and sequences $\langle C^\eta_j: j<\lambda\rangle$ for
$\eta\in {}^{\textstyle i}2$ such that 
\begin{enumerate}
\item[(a)] $\|B_\eta\|<\lambda$, $N_\eta\in\gK_\lambda$, $B_\eta\subseteq
  |N_\eta|\subseteq 
  \cup \{D_j : j < 1 + \lh (\eta ) \} $ 
\item[(b)] $\langle C^\eta_j:j<\lambda\rangle$ is increasing continuous,
  $\bigcup\limits_{j<\lambda} C^\eta_j=|N_\eta|$, $\|C^\eta_j\|<\lambda$,
\item[(c)] if $\nu\vartriangleleft\eta$ then $N_\nu\leq_\gK N_\eta$ and
  $B_\nu\subseteq B_\eta$,
\item[(d)] if $j_1,j_2<i$ then $C^{\eta\rest j_1}_{j_2}\subseteq B_\eta$,
\item[(e)] if the pair $(B_\eta,N_\eta)$ does not have the w.  amalgmation
  property in $\gK_\lambda$ then $N_{\eta\conc\langle 0\rangle}$, $N_{\eta
\conc\langle 1\rangle}$ witness it (i.e.\ they cannot be w.  amalgmated over
  $B_\eta$),
\item[(f)] if $i$ is limit and $\eta\in {}^{\textstyle i}2$ then $B_\eta=
\bigcup\limits_{j<i} B_{\eta\rest j}$, $\bigcup\limits_{j<i}N_{\eta\rest j}
\subseteq N_\eta$.
\end{enumerate}
There are no problems with carrying out the construction (remember that
$\leq_\gK$ is a nice partial order.
Finally, for $\eta\in {}^{\textstyle
\lambda}2$ we let $B_\eta=\bigcup\limits_{i<\lambda}B_{\eta\rest i}$ and
$N_\eta=\bigcup\limits_{i<\lambda} N_{\eta\rest i}$. Clearly, by
\ref{amalgdef}(2)$(\gamma$),
we have $N_\eta\in \gK$ and $B_\eta\subseteq
|N_\eta|$ for each $\eta\in {}^{\textstyle \lambda}2$. Moreover,
\[|N_\eta|=\bigcup_{j<\lambda}|N_{\eta\rest j}|=\bigcup_{j<\lambda}
\bigcup_{i<\lambda} C^{\eta\rest j}_i=\bigcup_{j^*<\lambda} \bigcup_{j_1,j_2<
j^*} C^{\eta\rest j_1}_{j_2}\subseteq\bigcup_{j^*<\lambda} B_{\eta\rest j^*}
=B_\eta,\] 
and thus $B_\eta=|N_\eta|$. Since $\gK$ is categorical in $\lambda$, for each
$\eta\in {}^{\textstyle \lambda} 2$ there is an isomorphism $f_\eta:N_\eta
\stackrel{\rm onto}{\longrightarrow} M$. 

Fix $\eta\in {}^{\textstyle\lambda}2$ for a moment.\\
Let $E_\eta=\{\delta<\lambda:f_\eta[B_{\eta\rest\delta}]=A_\delta=\delta\}$. 
Clearly, $E_\eta$ is a club of $\lambda$. Note that if $\delta\in E_\eta$
then:
\[\begin{array}{lcl}
(\boxtimes)\quad\delta\in S^{\bar{A}}_M&\Rightarrow&(A_\delta,M)\mbox{ does
  not have the w.  amalgmation property}\\
\ &\Rightarrow&(B_{\eta\rest\delta}, N_\eta)\mbox{ fails the w.  amalgmation
  property}\\ 
\ &\Rightarrow&(B_{\eta\rest\delta}, N_{\eta\rest\delta})\mbox{ fails the
  w.  amalgmation property}\\  
\ &\Rightarrow&N_{\eta\rest\delta\conc\langle 0\rangle},\ N_{\eta\rest\delta
  \conc\langle 1\rangle}\mbox{ cannot be w.  amalgmated}\\
\ &           &\qquad\qquad\mbox{ over }(B_{\eta\rest\delta},N_{\eta\rest
                \delta})\\   
\ &\Rightarrow&\mbox{for each }\nu\in{}^{\textstyle\lambda}2\mbox{ such that
  } \eta\rest\delta\conc\langle 1-\eta(\delta)\rangle\vartriangleleft\nu\\
\ &\ &\mbox{we have } f_\nu\rest B_{\eta\rest\delta}\neq f_\eta\rest B_{\eta
\rest\delta}.
\end{array}\]
We define a colouring 
\[F:\bigcup\limits_{\alpha<\lambda}{}^{\textstyle\alpha}(\Hla)\longrightarrow
\{0,1\}\]
by letting, for $f\in\DOM_\alpha$, $\alpha<\lambda$,
\[F(f)=1\quad\mbox{iff}\quad (\exists\eta\in{}^{\textstyle\lambda}2)\big(\eta(
\alpha)=0\ \&\ (\forall i<\alpha)(f(i)=(\eta(i),f^{-1}_\eta(i)))\big).\]
We have assumed $S^{\bar{A}}_M\notin\wdmid$, so there is $\rho\in{}^{
\textstyle \lambda}2$ such that for each $f\in\DOM_\lambda$ the set
\[
S_f=\{\delta\in S^{\bar{A}}_M:\rho(\delta)=F(f\rest\delta)\}
\]
is stationary and even $ \not=  \emptyset 
\mbox{ mod }
\idb(F^ \otimes )$.
Let $f\in\DOM_\lambda$ be defined by $f(i)=(\rho(i),f^{
-1}_\rho(i))$ (for $i<\lambda$). Note that

if $\alpha\in E_\rho$, $\rho(\alpha)=0$

then $\rho$ is a witness to $F(f\rest\alpha)=1$ and hence $\alpha\notin S_f$.

\noindent Since $S_f$ is stationary and even $ \not=  \emptyset 
\mod 
\idb(F^ \otimes ) $ 
and $E_\rho$ is a club of $\lambda$ we
may pick $\delta\in S_f\cap E_\rho$. Then $\rho(\delta)=1$ and hence
$F(f\rest\delta)=1$, so let $\eta_\delta\in {}^{\textstyle \lambda}2$ be a
witness for it. It follows from the definition of $F$ that then
$\eta_\delta(\delta)=0$, and $\eta_\delta\rest\delta=\rho\rest\delta$, and
$f^{-1}_{\eta_\delta}\rest\delta=f^{-1}_\rho\rest\delta$. 
    Hence
$f_{\eta_\eta }  
    \rest B_{\eta_\delta\rest\delta}= f_\rho\rest B_{\rho\rest
\delta}$, so both have range $A_\delta=\delta$ (and $\delta\in E_{
\eta_\delta}\cap E_\rho\cap S^{\bar{A}}_M$). But now we get a contradiction
with $(\boxtimes)$. 
\end{proof}

\bibliographystyle{amsalpha}
\bibliography{shlhetal}

\providecommand{\bysame}{\leavevmode\hbox to3em{\hrulefill}\thinspace}
\providecommand{\MR}{\relax\ifhmode\unskip\space\fi MR }
\providecommand{\MRhref}[2]{%
  \href{http://www.ams.org/mathscinet-getitem?mr=#1}{#2}
}
\providecommand{\href}[2]{#2}
\begin{thebibliography}{ALM13}

\bibitem[ALM13]{AsLaMo13}
David Aspero, Paul Larson, and Justin~Tatch Moore, \emph{{Forcing axioms and
  the continuum hypothesis}}, Acta Math \textbf{210} (2013), no.~1, 1--29.

\bibitem[DS78]{Sh:65}
Keith~J. Devlin and Saharon Shelah, \emph{{A weak version of $\diamondsuit$
  which follows from $2^{\aleph_0}<2^{\aleph_1}$}}, Israel J. Math. \textbf{29}
  (1978), no.~2-3, 239--247. \MR{0469756}

\bibitem[GS21]{Sh:1111}
Shimon Garti and Saharon Shelah, \emph{{Double weakness}}, Acta Math. Hungar.
  \textbf{163} (2021), no.~2, 379--391,
  \href{https://arxiv.org/abs/2002.03573}{arXiv: 2002.03573}. \MR{4227788}

\bibitem[Mak85]{Mw85}
Johann~A. Makowsky, \emph{Compactnes, embeddings and definability},
  Model-Theoretic Logics (J.~Barwise and S.~Feferman, eds.), Springer-Verlag,
  1985, pp.~645--716.

\bibitem[She77]{Sh:64}
Saharon Shelah, \emph{{Whitehead groups may be not free, even assuming CH. I}},
  Israel J. Math. \textbf{28} (1977), no.~3, 193--204. \MR{0469757}

\bibitem[She83]{Sh:87b}
\bysame, \emph{{Classification theory for nonelementary classes. I. The number
  of uncountable models of $\psi \in L_{\omega_1,\omega }$. Part B}}, Israel J.
  Math. \textbf{46} (1983), no.~4, 241--273. \MR{730343}

\bibitem[She85]{Sh:208}
\bysame, \emph{{More on the weak diamond}}, Ann. Pure Appl. Logic \textbf{28}
  (1985), no.~3, 315--318. \MR{790390}

\bibitem[She87a]{Sh:88}
\bysame, \emph{{Classification of nonelementary classes. II. Abstract
  elementary classes}}, {Classification theory (Chicago, IL, 1985)}, Lecture
  Notes in Math., vol. 1292, Springer, Berlin, 1987, pp.~419--497. \MR{1033034}

\bibitem[She87b]{Sh:192}
\bysame, \emph{{Uncountable groups have many nonconjugate subgroups}}, Ann.
  Pure Appl. Logic \textbf{36} (1987), no.~2, 153--206. \MR{911580}

\bibitem[She98]{Sh:f}
\bysame, \emph{{Proper and improper forcing}}, 2nd ed., Perspectives in
  Mathematical Logic, Springer-Verlag, Berlin, 1998. \MR{1623206}

\bibitem[She01]{Sh:576}
\bysame, \emph{{Categoricity of an abstract elementary class in two successive
  cardinals}}, Israel J. Math. \textbf{126} (2001), 29--128,
  \href{https://arxiv.org/abs/math/9805146}{arXiv: math/9805146}. \MR{1882033}

\bibitem[She08]{Sh:897}
\bysame, \emph{{Theories with Ehrenfeucht-Fra\"iss\'e equivalent non-isomorphic
  models}}, Tbil. Math. J. \textbf{1} (2008), 133--164,
  \href{https://arxiv.org/abs/math/0703477}{arXiv: math/0703477}. \MR{2563810}

\bibitem[She09a]{Sh:h}
\bysame, \emph{{Classification theory for abstract elementary classes}},
  Studies in Logic (London), vol.~18, College Publications, London, 2009.
  \MR{2643267}

\bibitem[She09b]{Sh:i}
\bysame, \emph{{Classification theory for abstract elementary classes. Vol.
  2}}, Studies in Logic (London), vol.~20, College Publications, London, 2009.
  \MR{2649290}

\bibitem[She09c]{Sh:88r}
\bysame, \emph{{Classification theory for elementary abstract classes}},
  Studies in Logic (London), vol.~18, College Publications, London, 2009,
  [Title on cover: Classification theory for abstract elementary classes],
  Mathematical Logic and Foundations
  \href{https://arxiv.org/abs/0705.4137}{arXiv: 0705.4137} Ch. I of [Sh:h].
  \MR{2643267}

\bibitem[She20]{Sh:1028}
\bysame, \emph{{Quite free complicated Abelian groups, pcf and black boxes}},
  Israel J. Math. \textbf{240} (2020), no.~1, 1--64,
  \href{https://arxiv.org/abs/1404.2775}{arXiv: 1404.2775}. \MR{4193126}

\bibitem[SZ99]{Sh:610}
Saharon Shelah and Jind{\v{r}}ich Zapletal, \emph{{Canonical models for
  $\aleph_1$-combinatorics}}, Ann. Pure Appl. Logic \textbf{98} (1999),
  no.~1-3, 217--259, \href{https://arxiv.org/abs/math/9806166}{arXiv:
  math/9806166}. \MR{1696852}

\end{thebibliography}

\end{document}